\documentclass[11pt,twoside,a4paper]{article}
\usepackage{bbm}
\usepackage{enumitem}
\usepackage{amsmath,mathtools,amsfonts,euscript,amssymb,amsthm,a4}
\usepackage{xcolor}
\usepackage[pagewise]{lineno}

\mathtoolsset{showonlyrefs}
\makeatletter \catcode`@=11
\newbox\tr@tto
\setbox\tr@tto=\hbox{{\count0=0\dimen0=-,9pt\dimen1=1,1pt\loop\ifnum
    \count0<11 \advance \count0 by1 \vrule width.51pt height\dimen1
    depth\dimen0\kern-0.17pt\advance\dimen0 by-0.05pt\advance\dimen1
    by0.1pt\repeat \loop\ifnum\count0<21\advance \count0 by1 \vrule
    width.6pt height\dimen1 depth\dimen0\kern-0.2pt \advance\dimen0
    by-0.1pt\advance\dimen1 by 0.05pt\repeat}}
\def\medint{\displaystyle\copy\tr@tto\kern-10.4pt\int}

 \allowdisplaybreaks

\def\R{\mathbb R}

\def\R{{\mathbb R}}
\def\RN{{\mathbb R}^{N}}
\def\M{{\mathbb R}^{n\times N}}

\def\Om{\Omega}
\def\N{\mathbb N}

\def\Phig{\mathbf{\Phi}}
\mathtoolsset{showonlyrefs}
\newcommand{\norm}[1]{\left\lVert#1\right\rVert}
\newtheorem{Thm}{Theorem}[section]
\newtheorem{Lem}[Thm]{Lemma}
\newtheorem{Rem}[Thm]{Remark}
\newtheorem{Prop}[Thm]{Proposition}
\newtheorem{Def}[Thm]{Definition}

\numberwithin{equation}{section}

\title{Quasiconvex Bulk 
and Surface Energies with subquadratic growth}

\author{
Menita Carozza\\
{\it  Department of Engineering, University of  Sannio}
\\ {\it Corso Garibaldi 107, 82100 Benevento, Italy}
\\ {\it e-mail: carozza@unisannio.it }
\bigskip
\\
  Luca Esposito\\
 {\it Dipartimento di Matematica, 
University of Salerno}\\ {\it Via Giovanni Paolo II, 
84084 Fisciano (SA),  Italy} \\{\it  e-mail: luesposi@unisa.it }
 \bigskip
 \\
  Lorenzo Lamberti\\
 {\it Institute \'Elie Cartan, Université de Lorraine, CNRS}\\ {\it F-54000 Nancy, France} \\{\it  e-mail: lorenzo.lamberti@univ-lorraine.fr}}

\begin{document}

\maketitle

\begin{abstract}
 \noindent 
We establish partial H\"older continuity of the gradient for equilibrium  configurations of {vectorial multidimensional} variational problems, involving bulk and surface energies. The bulk energy densities are uniformly strictly  quasiconvex functions  with $p$-growth, $1<p< 2$, without any further structure conditions. 
The anisotropic surface energy is defined by means of an elliptic integrand  $\Phi$ not necessarily regular.
\end{abstract}



\bigskip

\section{Introduction and statements}

Let us consider a functional $\mathcal{F}$ with density energy discontinuous through an interface $\partial A$, inside an open bounded subset $\Omega$  of $\R^n$, of the form
\begin{equation}\label{intro-1}
{\mathcal F}(v,A):=\int_\Om \left(F(D v)+\mathbbm{1}_{{A}}G(D v)\right)\,dx+P(A,\Om),
\end{equation}
where $v\in W^{1,p}_{loc}(\Om;\R^N)$,  and $F,G:\R^{n\times N}\rightarrow \R$ are  $C^2$-integrands.
Assume that these integrands satisfy the following growth and uniformly strict $p$-quasiconvexity conditions, for $p> 1$ and positive
constants $\ell_1,\,\ell_2,\,L_1,\,L_2$:

\begin{equation}
\label{F1p}\tag{$F1$}
    0\leq F(\xi)\leq L_1 (1+| \xi|^2 )^{\frac{p}{2}},
\end{equation}

\begin{equation}\label{F2p}\tag{$F2$}
    \int_{\Omega}F(\xi+D \varphi)\, dx \ge \int_{\Omega} \Big(F(\xi)+ \ell_1|D\varphi|^2(1+|D\varphi|^2)^{\frac{p-2}{2}}\Big)\,dx,
\end{equation}

\begin{equation}\label{G1p}\tag{$G1$}
    0\leq G(\xi)\leq  L_2 (1+| \xi|^2 )^{\frac{p}{2}},
\end{equation}

\begin{equation}\label{G2p}\tag{$G2$}
    \int_{\Omega} G(\xi+D \varphi)\, dx \ge \int_{\Omega} \Big(G(\xi)+ \ell_2|D \varphi|^2(1+|D\varphi|^2)^{\frac{p-2}{2}}\Big)\,dx,
\end{equation}
for every $\xi\in\M$ and $\varphi\in C_0^1(\Omega;\mathbb{R}^N)$.\\
\noindent Existence and regularity results have been obtained initially  in the scalar case ($N=1$) in \cite{AB,AR,CFP,E,EF,EL1,EL2,ELP,FJ,Lam,Lin,KL}. In the vectorial case ($N>1$), the authors in \cite{CFP1}  proved  the existence of local minimizers   of \eqref{intro-1}, for any $p>1$ under the quasiconvexity assumption quoted above. In the same paper,  the $C^{1,\alpha}$ partial regularity is proved for minimal configurations outside a negligible set, in the quadratic case $p=2$.\\
\noindent In  \cite {CEL} the same  regularity result has been established in the  general case $p\ge 2$, also addressing anisotropic surface energies. 
\noindent F.J. Almgren was the first to study such  surface energies  in his celebrated  paper \cite{A} (see also \cite{Bom,DS,F,SS,Simm} for subsequent results). This kind of energies arises in many physical contexts such as the formation of crystals (see \cite{BNP1,BNP2}), liquid drops (see \cite{CNT,FMA}), capillary surfaces (see \cite{DFM, DM}) and phase transitions (see \cite{Gur}).  \\
In this paper, we consider the same functional as in \cite {CEL}, given by
\begin{equation}\label{intro0}
{\mathcal I}(v,A):=\int_\Om \left(F(D v)+\mathbbm{1}_{_{A}}G(D v)\right)\,dx+\int_{\Omega\cap\partial^* A} \Phi(x,\nu_A(x))\;d\mathcal{H}^{n-1}(x),
\end{equation}
in the case of sub-quadratic growth,  $1<p<2$. We achieve analogous regularity results  as those established in \cite {CEL}, thereby completing the answer to the problem for all $p>1$.\\
\noindent In this setting $A\subset\Om$ is a  set of finite perimeter,  $u\in W_{\mathrm{loc}}^{1,p}(\Om;\mathbb{R}^N)$,  $\mathbbm{1}_A$ is the characteristic function of the set $A$,  $\partial^* A$ denotes the reduced boundary of $A$ in $\Omega$ and $\nu_A$ is the measure-theoretic outer unit
normal to $A$. Moreover,  $\Phi$ is an elliptic integrand on $\Omega$ (see Definition \ref{EI}), i.e. $\Phi:\overline{\Omega}\times \R^n\rightarrow [0,\infty]$  is lower semicontinuous, $\Phi(x,\cdot)$  is convex and positively one-homogeneous, $\Phi(x,t\nu)=t\Phi(x,\nu)$ for every $t\geq 0$,  and the  anisotropic surface energy of a set $A$ of finite perimeter in $\Omega$ is defined as follows
\begin{equation}\label{ASE0}
\Phig(A; B):=\int_{B\cap\partial^*A} \Phi(x,\nu_A(x))\;d\mathcal{H}^{n-1}(x),
\end{equation}
for every Borel set $B\subset\Omega$.
The further assumption 
\begin{equation}\label{nondege0}
\frac{1}{\Lambda}\leq \Phi(x,\nu)\leq \Lambda,
\end{equation}
with $\Lambda>1$, allows  to compare 
the surface energy introduced in \eqref{ASE0} with the usual perimeter.
\noindent Let us recall that in the vectorial setting, as in the previously cited papers, the regularity we can expect for the gradient of the minimal deformation $u:\Omega\rightarrow \R^N$, ($N> 1$), even in absence of a surface term, is limited to a partial regularity result.




\noindent We say that a pair $(u,E)$ is a local minimizer of ${\mathcal I}$ in $\Omega$, if for every open set $U\Subset \Omega$ and every pair $(v,A)$, where $v-u\in W_0^{1,p}(U;\mathbb{R}^N)$ and $A$ is a set of finite perimeter with $A\Delta E\Subset U$, we have
  \begin{equation*}\label{min}
  \int_U (F(Du)+\mathbbm{1}_{_{E}}G(D u))\,dx+\Phig(E;U)\leq \int_U (F(D v)+\mathbbm{1}_{_{A}}G(D v))\,dx+\Phig(A;U).
  \end{equation*}
 
\noindent Existence and  regularity results for local minimizers of integral functionals  with uniformly strict $p$-quasiconvex integrand, also in the non autonomous case, have been  widely investigated (see \cite{AF,AF1, cfm, cm, cp,cp2,GM,ma}  and  \cite{Gia,gi}).\\
\noindent Regarding the functional \eqref{intro0}, the existence of local minimizers   is guaranteed by the following theorem, proved in \cite {CEL}. 
\begin{Thm}\label{uno}
Let $p>1$ and assume that 
\eqref{F1p}, \eqref{F2p}, \eqref{G1p}, \eqref{G2p} hold.
Then, if $v\in W^{1,p}_{\mathrm{loc}}(\Omega;\mathbb{R}^N)$ and $A\subset \Om$ is a set of finite perimeter in $\Om$,  for every sequence $\{(v_k, A_k)\}_{k\in\N}$ such that $\{v_k\}$ weakly converges to $v$ in $W^{1,p}_{\mathrm{loc}}(\Omega;\mathbb{R}^N)$ and $\mathbbm{1}_{A_k}$ strongly converges to $\mathbbm{1}_{{A}}$ in $L^1_{\mathrm{loc}}(\Omega)$, we have
\begin{equation*}\label{lsc}{\mathcal I}(v,A)\leq \liminf_{k\to \infty}{\mathcal I}(v_k,A_k).\end{equation*}
In particular, ${\mathcal I}$ admits a minimal configuration $(u,\mathbbm{1}_E)\in W^{1,p}_{\mathrm{loc}}(\Omega;\mathbb{R}^N)\times BV_{\mathrm{loc}}(\Omega;[0,1])$
\end{Thm}
\noindent We emphasize that, in particular, the previous theorem implies the semicontinuity of the anisotropic perimeter functional \eqref{ASE0}.\\
\noindent In this paper, we obtain a   $C^{1,\alpha}$ regularity result for the minimizers of \eqref{intro0} in the case of sub-quadratic growth, $1<p<2$. If we further assume a closeness condition on $F$ and $G$ (see (H) in Theorem \ref{main}),   we prove that $u\in C^{1,\gamma}(\Om_1)$ for every   $\gamma\in (0,\frac{1}{p'})$ on a  full measure set $\Om_1\subset \Om$. Furthermore, we do not assume any regularity on  $\Phi$ in order to get the regularity of $u$. \\

 \noindent Our main theorem is the following,
\begin{Thm}\label{main}
 Let $(u,E)$ be a local minimizer of ${\mathcal I}$.
Let the bulk density energies $F$ and $G$ satisfy  $\eqref{F1p}$, $\eqref{F2p}$, $\eqref{G1p}$, $\eqref{G2p}$, with $1<p<2$, and  let the surface energy $\Phig$ be of general type \eqref{ASE0} with $\Phi$  satisfying \eqref{nondege0}.  Then there exist an exponent $\beta\in (0,1)$ and  an open set $\Omega_0\subset \Omega$ with full measure such that $u\in C^{1,\beta}(\Om_0;\mathbb{R}^N)$. If we assume in addition that
 \begin{equation}
 \label{H}\tag{H}
     \frac{L_2}{\ell_1+\ell_2}<1, 
 \end{equation}
  there exists an open set $\Omega_1\subset \Omega$ with full measure such that $u\in C^{1,\gamma}(\Om_1;\mathbb{R}^N)$ for every   $\gamma\in \big(0,\frac{1}{p'}\big)$.
\end{Thm}
\noindent 
The proof of the regularity of $u$ is based on a blow-up argument aimed to
establish a decay estimate for the excess function
\begin{equation}
 U(x_0,r):=\medint_{B_r(x_0)}\bigl|V(Du)- V\bigl((Du)_{x_0,r}\bigr)\bigr|^2\,dx + \frac{P(E, B_r(x_0))}{r^{n-1}}+r,
\end{equation}
where
\begin{equation*}
    V(\xi)=(1+|\xi|^2)^{(p-2)/4}\xi, \quad \forall\xi\in \R^k.
\end{equation*}
To this aim, we use a comparison argument between the blow-up sequence $v_h$ at small scale in the balls $B_{r_h}(x_h)$ and the  solution $v$ of a suitable linearized system. 
The challenging part of the argument, as usual, is to prove that the `good' decay estimates available for the function $v$ (see Proposition \ref{regu}), are inherited by the $v_h$ as $h\rightarrow \infty$.\\
To achieve this result,  the main tool is a Caccioppoli type inequality that we prove for minimizers of perturbed rescaled functionals (see \eqref{Cacciofin}) involving the function $V(Dv_h)$ and the perimeter of the rescaled minimal set $E_h$. The Caccioppoli inequality combined with the Sobolev Poincar\`e inequality will lead us to a contradiction (see Step 6 of Proposition \ref{decay1}). In this final step,  the issue to deal with the function $V(Du)$ in the sub-quadratic case, is overcome by using a suitable Sobolev Poincar\`e inequality involving  $V(Du)$ (see  Theorem \ref{SPsub}), whose proof is due to  \cite{cfm}.

\section{Preliminaries}

\bigskip

Let $\Omega$ be a bounded open set in $\R^n$,  $n\geq 2$ ,  $u:\Omega\rightarrow \R^N$, $N>1$. We denote by 
$
B_r(x):=\left\{y\in \R^n: |y-x|<r\right\}$  the open ball centered at $x\in \R^n$ of radius $r>0$, $\mathbb S^{n-1}$ represents the unit sphere of $\R^n$, $c$  a
generic constant that may vary.

\noindent For $B_r(x_0)\subset\R^n$ and $u\in L^1(B_r(x_0);\R^N)$ we denote 
\begin{equation*}
   (u)_{x_0,r}:= \medint_{B_r(x_0)}u(x)\,dx 
\end{equation*}
and we  will omit the dependence on the center when it is clear from the context.\\
\begin{equation*}
\langle \xi , \eta \rangle := \mbox{trace}(\xi^{T}\eta ),
\end{equation*}
for the usual inner product
of $\xi$ and $\eta$, and accordingly $| \xi | := \langle \xi , \xi \rangle^{\frac{1}{2}}$.
\noindent If $F \colon \M \to \R$ is sufficiently differentiable, we write
\begin{equation*}
    D F(\xi )\eta :=\sum_{\alpha=1}^N\sum_{i=1}^n\frac{\partial F}{\partial \xi^\alpha_i}(\xi)\eta^\alpha_i
\quad \text{ and } \quad
D^2 F(\xi )\eta\eta := \sum_{\alpha,\beta=1}^N\sum_{i,j=1}^n\frac{\partial F}{\partial \xi_i^\alpha \partial \xi_j^\beta}(\xi)\eta^\alpha_i\eta^\beta_j,
\end{equation*}
for $\xi$, $\eta \in \M$.\\
\noindent  It is well known that for quasiconvex $C^1$ integrands the assumptions \eqref{F1p} and \eqref{G1p} yield the upper bounds

\begin{equation}\label{(H4)}
|D_\xi F(\xi)|\leq c_1L_1(1+|\xi|^2)^{\frac{p-1}{2}}
\quad\mathrm{and}\quad
|D_\xi G(\xi)|\leq  c_2L_2(1+|\xi|^2)^{\frac{p-1}{2}}
\end{equation}
for all $\xi \in \M$, with $c_1 $ and $c_2$  constants depending only on $p$ (see \cite[Lemma 5.2]{gi} or \cite{ma}).\\
\noindent Furthermore, if $F$ and $G$ are $C^2$, then
$\eqref{F2p}$ and $\eqref{G2p}$ imply  the following  strong Legendre-Hadamard conditions
\begin{equation*}
 \sum_{\alpha,\beta=1}^N\sum_{i,j=1}^n\frac{\partial F}{\partial \xi_i^\alpha \partial \xi_j^\beta}(Q) \lambda_i\lambda_j\mu^\alpha\mu^\beta\ge c_3|\lambda|^2|\mu|^2
\quad\text{and}\!\quad\!
\sum_{\alpha,\beta=1}^N\sum_{i,j=1}^n\frac{\partial G}{\partial \xi_i^\alpha \partial \xi_j^\beta}(Q) \lambda_i\lambda_j\mu^\alpha\mu^\beta\ge  c_4|\lambda|^2|\mu|^2, 
\end{equation*}
for all $Q\in\M$, $\lambda\in \R^n$, $\mu \in \R^N $, where $c_3=c_3(p,\ell_1)$ and $c_4=c_4(p,\ell_2)$ are positive constants (see \cite[Proposition 5.2]{gi}).

We will  need the following quite standard regularity result (see \cite{cfm} for its proof).

\begin{Prop}\label{regu}
Let $v\in W^{1,1}(\Omega;\R^N)$ be such that
\begin{equation*}
    \int_\Omega Q^{ij}_{\alpha \beta} D_i v^\alpha D_j \varphi^\beta\,dx=0,
\end{equation*}
for every $\varphi \in C_{c}^\infty(\Omega;\R^N)$, where $Q=\{Q^{ij}_{\alpha\beta}\}$
is a constant matrix satisfying $|Q^{ij}_{\alpha \beta}|\leq L$ and the strong Legendre-Hadamard condition
\begin{equation*}
   Q^{ij}_{\alpha\beta}\lambda_i\lambda_j\mu^\alpha\mu^\beta\geq\ell|\lambda|^2|\mu|^2, 
\end{equation*}
for all $\lambda\in\R^n$, $\mu\in\R^N$ and for some positive constants $\ell,L>0$.
Then $v \in C^\infty$ and, for any $B_ R (x_0 ) \subset \Omega$, the following estimate holds
\begin{equation}\label{supest}
\sup_{B_{R/2}}|D v|\leq \frac{c}{R^n}\int_{B_R}|D v|\,dx,
\end{equation}
where $c=c(n,N,\ell,L)>0$ .
\end{Prop}
\noindent 
We assume that $1<p<2$ and we  refer to the auxiliary function 
\begin{equation}\label{V}
    V(\xi)=(1+|\xi|^2)^{(p-2)/4}\xi, \quad \forall\xi\in \R^k,
\end{equation}
whose useful properties are listed in the following lemma (see \cite{cfm} for the proof).
\begin{Lem}\label{PropV} Let $1<p<2$ and let $V:\R^k\rightarrow \R^k$ be the function defined in \eqref{V}, then for any $\xi,\eta \in \R^k$ and $t>0$ the following inequalities hold:
\begin{itemize}
    \item[(i)]\; $2^{(p-2)/4}\min\{|\xi|,|\xi|^{p/2}\}\leq \bigl|V(\xi)\bigr|\leq \min\{|\xi|,|\xi|^{p/2}\}$,
    \item[(ii)]\; $\bigl|V(t\xi)\bigr|\leq \max\{t,t^{p/2}\}\bigl|V(\xi)\bigr|$ ,
    \item[(iii)]\; $\bigl|V(\xi+\eta)\bigr|\leq c(p)\bigl[\bigl|V(\xi)\bigr|+\bigl|V(\eta)\bigr|\bigr]$,
    \item[(iv)]\; $\frac p2|\xi-\eta|\leq \bigl(1+|\xi|^2+|\eta|^2\bigr)^{(2-p)/4}\bigl|V(\xi)-V(\eta)\bigr|\leq c(k,p)|\xi-\eta|$,
    \item[(v)]\; $\bigl|V(\xi)-V(\eta)\bigr|\leq c(k,p)\bigl|V(\xi-\eta)\bigr|$,
    \item[(vi)]\; $\bigl|V(\xi-\eta)\bigr|\leq c(p,M)\bigl|V(\xi)-V(\eta)\bigr|$, if $|\eta| \leq M$.
\end{itemize}
\end{Lem}

\medskip

We will also use the following iteration lemma  (see  \cite[Lemma 6.1]{gi})
\begin{Lem}\label{widman}
Let $0<\rho<R $ and let $\psi \colon [\rho,R] \to \R$ be a bounded non negative function. Assume that for all $\rho \leq s < t \leq R$
we have
\begin{equation}
\psi (s) \leq \vartheta \psi (t) + A + \frac{B}{(s-t)^\alpha}+\frac{C}{(s-t)^\beta}
\end{equation}
where $\vartheta \in [0,1)$, $\alpha>\beta>0$ and $A,B,C \geq 0$  are constants.
Then there exists a constant $c= c(\vartheta,\alpha)>0$ such that
\begin{equation}
\psi \bigl( \rho \bigr) \leq
c\left( A + \frac{B}{(R-\rho)^\alpha}+\frac{C}{(R-\rho)^\beta} \right).
\end{equation}
\end{Lem}
\noindent
 An easy extension of this result can be obtained by replacing homogeneity with condition (ii) of Lemma \ref{PropV}.
\begin{Lem}\label{iterationV}
Let $R>0$ and let $\psi \colon [R/2,R] \to [0,+\infty)$ be a bounded function. Assume that for all $R/2 \leq s < t \leq R$
we have
\begin{equation*}
\psi (s) \leq \vartheta \psi (t) + A \int_{B_R}\bigg|V\bigg(\frac{h(x)}{t-s}\bigg)\bigg|^2 dx +B, 
\end{equation*}
where $h\in L^p(B_r)$, $A,B>0$, and $0<\vartheta<1$. Then there exists a constant $c(\vartheta)>0$ such that
\begin{equation*}
\psi \bigg( \frac{R}{2} \bigg) \leq
 c(\vartheta)\bigg( A \int_{B_R}\Bigl|V\bigg(\frac{h(x)}{R}\bigg)\Bigr|^2 dx +B\bigg).
\end{equation*}
\end{Lem}

\medskip
\noindent Given a $C^1$ function $f:\R^k\to  \R$,   $Q\in \R^k$ and $\lambda>0$, we set
\begin{equation}\label{dfg}
f_{Q,\lambda}(\xi):=\frac{f(Q+\lambda\xi)-f(Q)-D f(Q)\lambda\xi}{\lambda^2},\quad\forall\xi\in \R^k.
\end{equation}
In the next sections we will use the following lemma about the growth of $f_{Q,\lambda}$ and $Df_{Q,\lambda}$.
\begin{Lem}\label{nic}
   Let $1<p<\infty$, and let $f$ be a $C^{2}(\R^k)$ function such that \begin{equation*}
       | f(\xi)| \leq  L \big(1+| \xi|^p)\quad \text{and}\quad
 |D f(\xi)|  \leq  L \big(1+|\xi |^2\big)^{(p-1)/2},
   \end{equation*}
 for any $\xi\in\R^k$ and for some $L>0$. Then for every $M>0$ there exists a constant $c=c(p,L,M)>0$ such that, for every $Q\in \R^k$, $|Q|\leq M$ and $\lambda>0$, it holds
\begin{equation}\label{R1}
 | f_{Q,\lambda}(\xi)| \leq  c \big(1+|\lambda \xi|^2\big)^{(p-2)/2}|\xi|^2\quad and\quad
 |D f_{Q,\lambda}(\xi)|  \leq  c\big(1+|\lambda \xi|^2\big)^{(p-2)/2}|\xi|,
\end{equation}
for all $\xi\in\R^k$.
\end{Lem}
\begin{proof}
Applying Taylor's formula, there exists $\theta\in [0,1]$ such that, for every $\xi \in \R^k$,
\begin{equation}
  f_{Q,\lambda}(\xi)=\frac 12 D^2 f(Q+\theta\lambda \xi)\xi\xi,  
\end{equation}
\begin{equation}
    Df_{Q,\lambda}(\xi)=\frac{1}{\lambda}\bigl(D f(Q+\lambda\xi)-D f(Q)\bigr)=\int_0^1 D^2 f(Q+\theta\lambda \xi)\xi \, d\theta.
\end{equation}
If we denote $K_M:=\max\left\{|D^2 f(\xi)|: |\xi|\leq M+1\right\}$, we have 
\begin{equation}\label{Co1}
|f_{Q,\lambda}(\xi)|\leq \frac 12 K_M |\xi|^2,\quad|D f_{Q,\lambda}(\xi)|\leq K_M |\xi|,\quad \text{if } |\lambda\xi|\leq 1.
\end{equation}
On the other hand, using growth condition \eqref{R1} and the definitions of $f_{Q,\lambda}$ and $Df_{Q,\lambda}$, we get 
\begin{equation}\label{Co2}
|f_{Q,\lambda}(\xi)|\leq c(p,L,M)\lambda^{p-2}|\xi|^p,\quad|D f_{Q,\lambda}(\xi)|\leq c(L,M)\lambda^{p-2}|\xi|^{p-1},\quad \text{whereas }|\lambda\xi|> 1.
\end{equation}
We get the result by combining \eqref{Co1} and \eqref{Co2} .
\end{proof}
\noindent A fundamental tool in order to handle the subquadratic  case is the following Sobolev-Poincar\'e inequality related to the function $V$, proved in  \cite{cfm}.
\begin{Thm} \label{SPsub}
If $1<p<2$, there exist $2/p <\alpha< 2$ and $\sigma > 0$ such that if $u \in W^{1,p}(B_{3R}(x_0), \R^N)$, then
\begin{equation}\label{SoPoV}
    \Bigg( \medint_{B_{R}(x_0)} \bigg|V\bigg(\frac{u - u_{x_o,R}}{R}\bigg)\bigg|^{2(1+\sigma)} dx \Bigg)^{\frac{1}{2(1+\sigma)}} 
    \leq C \bigg( \medint_{B_{3R}(x_0)} \bigl|V\bigl(Du\bigr)\bigr|^{\alpha} dx \bigg)^{\frac{1}{\alpha}},
\end{equation}
where the positive constant $C = C(n, N, p)$ is independent of $R$ and $u$.
\end{Thm}

\subsection{Sets of finite perimeter and anisotropic surface energies}
In this subsection we recall some elementary definitions and well-known properties of sets of finite perimeter. We introduce the notion of anisotropic perimeter as well.\\
\noindent Given a set $E\subset \R^n$ and $t\in [0,1]$, we define the set of points of $E$ of density $t$ as
\begin{equation}
E^{(t)}=\left\{x\in \R^n:\ |E\cap B_r(x)|=t|B_r(x)|+o(r^n)\text{ as }r\rightarrow 0^+\right\}.
\end{equation}
Let $U$ be an open subset $U$ of $\R^n$. A Lebesgue measurable set $E\subset \R^n$
is said to be a set of locally finite perimeter in $U$ if there exists a $\R^n$-valued Radon measure $\mu_E$ on $U$ (called the Gauss-Green measure of $E$) such that
\begin{equation*}
\int_{E}\nabla \phi\ dx=\int_{U}\phi \, d\mu_E,\quad\forall \phi \in C^1_c(U).
\end{equation*}
Moreover, we denote the perimeter of $E$ relative to $G\subset U$ by $P(E,G)=|\mu_E|(G)$.\\
 \noindent It is well known that the support of $\mu_E$ can be characterized by
\begin{equation}\label{support}
\text{spt}\mu_E=\bigl\{x\in U: 0<|E\cap B_r(x)|<\omega_n r^n, \,\forall r>0\bigr\}\subset U\cap \partial E,
\end{equation}
(see \cite[Proposition 12.19]{Ma}).
If $E$ is of finite perimeter in $U$, the {\it reduced boundary} $\partial^*E \subset U$ of $E$ is the set of those $x\in U$ such that
\begin{equation}\label{RB}
\nu_E(x):=\lim_{r\rightarrow 0^+}\frac{\mu_E(B_r(x))}{|\mu_E|(B_r(x))}
\end{equation}
exists and belongs to $\mathbb S^{n-1}$. The {\it essential boundary} of $E$ is defined as $\partial^e E:=\R^n\setminus(E^0\cup E^1)$. It is well-understood that
\begin{equation*}
\partial^* E\subset U\cap\partial^e E\subset {\text{spt}}\mu_E\subset U\cap\partial E,\hspace{1cm}U\cap\overline{\partial^*E}={\text{spt}}\mu_E.
\end{equation*}
Furthermore, Federer's criterion (see for instance \cite[Theorem 16.2]{Ma}) ensures that
\begin{equation}\label{Federer}
\mathcal{H}^{n-1}((U\cap\partial^e E)\setminus \partial^*E)=0.
\end{equation}
By De Giorgi’s rectifiability theorem (see \cite[Theorem 15.9]{Ma}), the Gauss-Green measure $\mu_E$ is completely characterized as follows:
\begin{equation}\label{DGRep}
\mu_E=\nu_E\mathcal{H}^{n-1}\llcorner\partial^*E, \quad |\mu_E|=\mathcal{H}^{n-1}\llcorner\partial^*E.
\end{equation}
The equality holds in the class of Borel sets compactly contained in $U$. Here, we have denoted $\mu\llcorner \partial^*E(F) = \mu (\partial^*E \cap F)$, for any subset $F$ of $\R^n.$\\
\begin{Rem}[Minimal topological boundary]\label{MTB}
If $E\subset\R^n$ is a set of locally finite perimeter in $U$ and $F\subset\R^n$ is such that $|(E\Delta F)\cap U|=0$, then $F$ is a set of locally finite perimeter in $U$ and $\mu_E=\mu_F$. In the rest of the paper, the
topological boundary $\partial E$ must be understood by considering the suitable representative of $E$ in order to have that $\overline{\partial^*E}=\partial E\cap U$. We will choose $E^{(1)}$ as representative of $E$. With such a choice it can be easily verified that 
\begin{equation}
U\cap\partial E=\bigl\{x\in U: 0<|E\cap B_r(x)|<\omega_n r^n, \forall r>0\bigr\}.
\end{equation}
Therefore, by $\eqref{support}$,
\begin{equation*}
\overline{\partial^*E}={\textnormal{spt}}\mu_E=\partial E\cap U.
\end{equation*}
\end{Rem}
\noindent In what follows, we give the definition of anisotropic surface energies and we recall some properties.
\begin{Def}[Elliptic integrands]\label{EI}
Given an open subset $\Omega$ of $\R^n$, $\Phi:\overline{\Omega}\times \R^n\rightarrow [0,\infty]$ is said to be an elliptic integrand on $\Omega$ if it is lower semicontinuous, with $\Phi(x,\cdot)$ convex and positively one-homogeneous for any $x\in\overline{\Omega}$, i.e. $\Phi(x,t\nu)=t\Phi(x,\nu)$ for every $t\geq 0$. Accordingly, the anisotropic surface energy of a set $E$ of finite perimeter in $\Omega$ is defined as
\begin{equation}\label{ASE}
\Phig(E;B):=\int_{B\cap\partial^*E} \Phi(x,\nu_E(x))\;d\mathcal{H}^{n-1}(x),
\end{equation}
for every Borel set $B\subset\Omega$.
\end{Def}
\noindent In order to prove the regularity of minimizers of anisotropic surface energies, it is well known that a $C^k$-dependence of the integrand $\Phi$ on the variable $\nu$, and a continuity condition with respect to the variable $x$, must be assumed (see the seminal paper \cite{A}). In fact, one more condition is essential, that is a non-degeneracy type condition for the integrand $\Phi$. More precisely, we have to assume that there exists a constant $ \Lambda>1$ such that
\begin{equation}\label{nondege}
\frac{1}{\Lambda}\leq \Phi(x,\nu)\leq \Lambda,
\end{equation}
for any $x\in\Omega$ and $\nu\in\mathbb{S}^{n-1}$.
We emphasize that \eqref{nondege} is the only assumption we make for the elliptic integrand $\Phi$. We observe that, if the elliptic integrand $\Phi$ satisfies the previous condition, then the anisotropic surface energy \eqref{ASE} satisfies the following comparability condition to the perimeter:
\begin{equation}
\frac{1}{\Lambda}\mathcal{H}^{n-1}(B\cap\partial^*E)\leq\Phig(E;B)\leq \Lambda \mathcal{H}^{n-1}(B\cap\partial^*E),
\end{equation}
for any set $E$ of finite perimeter in $\Omega$ and any Borel set $B\subset\Omega$.
\noindent A  useful relation is given by  proposition below proved in \cite{CEL}.
\begin{Prop}
\label{PerimetroAnisotropoUnione}
Let $U\subset\R^n$ be an open set and let $E,F\subset U$ be two sets of finite perimeter in $U$. It holds that
\begin{equation}\label{AdP}
\Phig(E\cup F;U)=\Phig(E;F^{(0)})+\Phig(F;E^{(0)})+\Phig(E;\left\{\nu_E=\nu_F\right\}).
\end{equation}
\end{Prop}

\medskip
\bigskip

\section{Decay Estimates}

\bigskip

In this section we prove decay estimates for minimizers of  functionals  \eqref{intro0} by using a well-known blow-up technique involving a suitable excess function. We consider the 
bulk
excess function defined as
\begin{equation}\label{excess}
U(x_0,r):=\medint_{B_r(x_0)}\bigl|V(Du)- V\bigl((Du)_{x_0,r}\bigr)\bigr|^2\,dx,
     \end{equation}
for $B_r(x_0)\subset\Omega$.
 
 \noindent When the assumption $\eqref{H}$ is in force, we refer to the following $\lq\lq$hybrid" excess
\begin{equation}\label{excess1}
U_*(x_0,r):=
U(x_0,r)
+ \frac{P(E, B_r(x_0))}{r^{n-1}}+r.
\end{equation}

\medskip
\begin{Prop}\label{decay1}  Let $(u,E)$ be a local minimizer of 
the functional $\mathcal{I}$  in \eqref{intro0} and let the assumptions $\eqref{F1p}$, $\eqref{F2p}$, $\eqref{G1p}$, $\eqref{G2p}$ and \eqref{H} hold. For every $M>0$ and every $0<\tau <\frac{1}{4}$, there exist two constants $\varepsilon_0=\varepsilon_0(\tau,M)>0$ and $c_*=c_*(n,p,\ell_1,\ell_2,L_1,L_2,\Lambda,M)>0$ such that if in $B_r(x_0)\Subset\Om$ it holds 
\begin{equation*}
     |(Du)_{x_0,r}|\leq M\quad\mathrm{and}\quad  U_*(x_{0}, r)\leq \varepsilon_0,
\end{equation*}
then
\begin{equation}\label{decay3}
U_*(x_{0}, \tau r)\leq c_* \tau U_*(x_{0}, r).
\end{equation}
\end{Prop}

\begin{proof}
In order to prove \eqref{decay3}, we argue by contradiction. Let $M>0$ and $ \tau\in (0,1/4)$ be such that  for every $h\in \mathbb{N}$, $C_*>0$, there exists a   ball  $ B_{ r_h}(x_h)\Subset \Omega$ such that
\begin{equation}\label{decay04}|(Du)_{x_h,r_h}|\leq M,\quad U_*(x_h,r_h)\to 0\end{equation}
and
\begin{equation}\label{decay4}
U_*(x_h,\tau r_h)\ge C_*\tau U_*(x_h,r_h).
\end{equation}
The constant $C_*$  will be determined later. We remark that we can confine ourselves to the case in which $E\cap  B_{ r_h}(x_h)\not=\emptyset$, since the case in which  $ B_{ r_h}(x_h)\subset \Omega\setminus E$ is easier, being $U=U_*$.

\medskip

\noindent{\bf Step 1.}\,\,{\it  Blow-up.}

\medskip

\noindent  We set $\lambda_h^2:=U_*(x_h,r_h)$,
$A_h:=(Du)_{x_h,r_h}$, $a_h:=(u)_{x_h,r_h}$, and we define
\begin{equation}\label{23}
v_h(y):=\frac{u(x_h+ r_hy)-a_h-r_hA_h y}{\lambda_h r_h}, \quad \forall y\in B_1.
\end{equation}
One can easily check that $(Dv_h)_{0,1}=0$
and $(v_h)_{0,1}=0$.
\noindent We set
\begin{equation*}
    E_h:= \frac{E-x_h}{r_h},\quad E^*_h:= \frac{E-x_h}{r_h}\cap B_1.
\end{equation*}
 By
 using $(ii)$ and $(vi)$ of Lemma \ref{PropV}, we deduce
\begin{align*}
\medint_{B_1}|V(Dv_h(y))|^2\,dy\leq &
\medint_{B_{r_h}(x_h)}\bigg|V\bigg(\frac{Du(x)- (Du)_{x_h,r_h}}{\lambda_h}\bigg)\bigg|^2\,dx \\
\leq &\frac{c(M)}{\lambda_h^2}
\medint_{B_{r_h}(x_h)}\bigl|V(Du(x))- V\bigl((Du)_{x_h,r_h}\bigr)\bigr|^2\,dx.
\end{align*}
Then, since
\begin{align}\label{resc}
\lambda_h^2 = U_*(x_h,r_h)
& =\medint_{B_1}\bigl|V\bigl(Du(x_h+r_hy)\bigr)-V\bigl( A_h\bigr)\bigr|^2\,dy+ \frac{P(E, B_{r_h}(x_h))}{r_h^{n-1}}+r_h,
\end{align}
it follows that $r_h\to 0$, $P(E_h,B_1)\to 0,$ and
\begin{equation}\label{24bis}
\frac{r_h}{\lambda_h^2}\leq 1,\quad \medint_{B_1}\bigl|V\bigl(Dv_h(y)\bigr)\bigr|^2\,dy
\leq c(M), \quad \frac{P(E_h,B_1)}{\lambda_h^2}
\leq 1.
\end{equation}
Therefore, by  \eqref{decay04}  and \eqref{24bis},  there exist a (not relabeled) subsequence of $\{v_h\}_{h\in\N}$, $A\in \M$ and $v\in W^{1,p}(B_1;\RN)$, such that
\begin{align}\label{25}
 & v_h\rightharpoonup v \quad\text{weakly in }W^{1,p}(B_{1}; \mathbb{R}^{N}),\quad v_h\to v\quad \text{strongly in }L^{p}(B_{1}; \mathbb{R}^{N}),\\
 & A_h\to A,\quad \lambda_h Dv_h\rightarrow 0 \quad\text{in }L^{p}(B_{1};\R^{n\times N}) \text{ and pointwise a.e. in }B_1,
 \end{align}
where we have used the fact that $(v_h)_{0,1}=0$.
Moreover, by  \eqref{24bis} and \eqref{decay04}, we have 
\begin{equation}\label{25quater}
\lim_{h\rightarrow\infty}\frac{(P(E_h,B_1))^{\frac{n}{n-1}}}{\lambda_h^2}\leq \lim_{h\rightarrow\infty} (P(E_h,B_1))^{\frac{1}{n-1}}\limsup_{h\rightarrow\infty}\frac{P(E_h,B_1)}{\lambda_h^2}=0.
\end{equation}
Therefore, by the relative isoperimetric inequality,
\begin{equation}\label{25ter}\lim_{h\rightarrow\infty}\min\left\{\frac{|E^*_h|}{\lambda_h^2}, \frac{|B_1\setminus E_h|}{\lambda_h^2}\right\}\leq c(n) \lim_{h\rightarrow\infty}\frac{\left(P(E_h,B_1)\right)^{\frac{n}{n-1}}}{\lambda_h^2}=0.
\end{equation}

\noindent In the sequel the proof will proceed differently depending on  
\begin{equation*}
\min\{|E^*_h|, |B_1\setminus E_h|\}=|E^*_h| \hspace{10mm}\text{ or } \hspace{10mm}
\min\{|E^*_h|, |B_1\setminus E_h|\}=|B_1\setminus E_h|.
\end{equation*}
The first case is easier to handle. To understand the reason, let us introduce the  expansions of  $F$ and $G$ around $A_h$ as follows:
\begin{align}\label{risc}
& F_h(\xi)
:=\frac{F(A_h+\lambda_h\xi)-F(A_h)-D
F(A_h)\lambda_h \xi}{\lambda_h^2},\\
& G_h(\xi):= \frac{G(A_h+\lambda_h\xi)-G(A_h)-D
G(A_h)\lambda_h \xi}{\lambda_h^2},
\end{align}
for any $\xi\in\R^{n\times N}$. In the first case the suitable rescaled functional to consider in the blow-up procedure is the following

\begin{equation}\label{26}
{\cal I}_h(w):=\int_{B_1}\big[F_h(Dw)dy +\mathbbm{1}_{E_h^*}G_h(Dw)\big]\,dy.
\end{equation}
We claim that $v_h$ satisfies the minimality inequality
\begin{equation}\label{29}
{\cal I}_h(v_h) \leq {\cal
I}_h(v_h+\psi)+\frac{1}{\lambda_h}\int_{B_1}\mathbbm{1}_{E_h^*} D G(A_h)D\psi(y)\,dy,
\end{equation}
for any $\psi\in W^{1,p}_0(B_1;\R^N)$. Indeed, using  the minimality of $(u,E)$ with respect to $(u+\varphi,E)$, for $\varphi\in W^{1,p}_0(B_{r_h}(x_h);\R^N)$, the change of variable $x=x_h+ r_hy$, setting $\psi(y):=\frac{\varphi(x_h+ r_hy)}{\lambda_h r_h}$, it holds that
\begin{align}\label{decay5a}
& \int_{B_1}\big[(F_h(D v_h(y))+\mathbbm{1}_{E_h^*}G_h(D v_h(y))\big]\,dy\notag\\
& \leq \int_{B_1}\big[F_h(D v_h(y)+D \psi(y))+\mathbbm{1}_{E_h^*}G_h(D v_h(y)+D\psi(y))\big]\,dy+\frac{1}{\lambda_h}\int_{B_1}\mathbbm{1}_{E_h^*} D G(A_h)D\psi(y)\,dy,
\end{align}
and \eqref{29} follows by the definition of ${\cal I}_h$ in \eqref{26}.\\
\noindent In the second case, the suitable rescaled functional to consider in the blow-up procedure is 
\begin{equation}\label{27}
{\cal H}_h(w):=\int_{B_1}\big[F_h(Dw) +G_h(Dw)\big]\,dy.
\end{equation}
We claim that
\begin{align}\label{decay5a00}
{\cal H}_h(v_h)\leq {\cal H}_h(v_h+\psi)+\frac{L_2}{\lambda_h^2}\int_{(B_1\setminus E_h)\cap \mathrm{supp}\psi}(\mu^2+|A_h+\lambda_hDv_h|^2)^{\frac{p}{2}}\,dy,
\end{align}
for all $\psi\in W^{1,p}_0(B_1;\R^{N})$.
Indeed,  the minimality of $(u,E)$ with respect to $(u+\varphi,E)$, for $\varphi\in W^{1,p}_0(B_{r_h}(x_h);\R^{N})$, implies that
\begin{align}\label{minHa}
& \int_ {B_{r_h}(x_h)} (F+G)(Du)\,dx=\int_{B_{r_h}(x_h)} \big[F(Du)+\mathbbm{1}_E G(Du)\big]\,dx+\int_{B_{r_h}(x_h)\setminus E}G(Du)dx\notag\\
& \leq \int_{B_{r_h}(x_h)}\big[ F(Du+D\varphi)+\mathbbm{1}_E G(Du+D\varphi)\big]\,dx+\int_{B_{r_h}(x_h)\setminus E}G(Du)dx\notag\\
& =\int_{B_{r_h}(x_h)} (F+G)(Du+D\varphi)dx+\int_{B_{r_h}(x_h)\setminus E}\big[G(Du)-G(Du+D\varphi)\big]\,dx\notag\\
& \leq \int_{B_{r_h}(x_h)} (F+G)(Du+D\varphi)dx +\int_{(B_{r_h}(x_h)\setminus E)\cap \mathrm{supp}\varphi}G(Du)dx ,
\end{align}
where we used that the last integral vanishes outside the support of $\varphi$ and that $G\ge 0$.
Using  the change of variable $x=x_h+ r_hy$ in the previous formula, we get
\begin{align*}
\int_{B_1}(F+G)(Du(x_h+r_hy))dy
& \leq \int_{B_1}(F+G)(Du(x_h+r_hy)+D\varphi(x_h+r_hy))\,dy\\
& +\int_{(B_1\setminus E_h)\cap \mathrm{supp}\psi} G(Du(x_h+r_hy))dy,
\end{align*}
or, equivalently, using  the definitions of $v_h$,
\begin{align*}\int_{B_1}(F+G)(A_h+\lambda_hDv_h)dy
& \leq\int_{B_1}(F+G)(A_h+\lambda_h(Dv_h+D{\psi}))\,dy\notag\\
& +\int_{(B_1\setminus E_h)\cap \mathrm{supp}\psi}G(A_h+\lambda_hDv_h)dy
\end{align*}
where ${\psi(y)}:=\frac{\varphi(x_h+r_hy)}{\lambda_h r_h}$, for $y\in B_1$.
Therefore, setting  
\begin{equation*}
    H_h:=F_h+G_h,
\end{equation*}
by the definitions of $F_h$ and $G_h$ in \eqref{risc} and using the assumption $\eqref{G1p}$, we have that
\begin{align}\label{decay5aaa} \int_{B_1} H_h(Dv_h)dy
& \leq\int_{B_1} H_h(Dv_h+D\psi)dy
+\frac{1}{\lambda_h^2}\int_{(B_1\setminus E_h)\cap \mathrm{supp}\psi}\!\!\! G(A_h+\lambda_hDv_h)\,dy\notag\\
& \leq\int_{B_1}H_h(Dv_h+D\psi)\,dy+\frac{L_2}{\lambda_h^2}\int_{(B_1\setminus E_h)\cap \mathrm{supp}{\psi}}\!\!\!\!\!\!\!\!\!\!\!\!\!\!\!\!\!\!\!\!\big(1+|A_h+\lambda_hDv_h|^2\big)^{\frac{p}{2}}\,dy,
\end{align}
i.e. \eqref{decay5a00}.
\medskip

\noindent{\bf Step 2.}\,\,{\it  A Caccioppoli type inequality.}

\medskip

 \noindent We claim that there exists a constant $c=c(n,p,\ell_1,\ell_2,L_1,L_2,M)>0$ such that for every $0<\rho<1$ there exists $h_0=h_0(n,p,M,\rho)\in \mathbb{N}$ such that 
\begin{align}\label{Cacciofin}
& \int_{B_{\frac{\rho}{2}}}\big|V\big(\lambda_h(Dv_h -(Dv_h)_{\frac{\rho}{2}}\big)\big|^2\,dy\\
& \leq c\Bigg[\int_{B_\rho}\bigg|V\bigg(\frac{\lambda_h\big(v_h -(v_h)_\rho-(Dv_h)_{\frac{\rho}{2}}\,y\big)}{\rho}\bigg)\bigg|^2\,dy+P(E_h,B_1)^{\frac{n}{n-1}}\Bigg]\notag,
\end{align}
for all $h>h_0$.
\noindent We  divide the proof  into  two  steps.

\medskip

\noindent{\bf Substep 2.a}\,\,{\it  The case $\min\{|E^*_h|, |B_1\setminus E_h|\}=|E^*_h|$.} 
\medskip

\noindent We consider $0<\frac{\rho}{2}<s<t<\rho<1$ and let $\eta\in C^\infty_0(B_{t})$  be a cut off function between $B_s$ and $B_t$, i.e. $0\leq \eta\leq 1$, $\eta\equiv 1$ on $B_s$ and $|\nabla \eta |\leq \frac{c}{t-s}$.
Set $b_h := (v_h)_{B_\rho}$ , $B_h := (Dv_h)_{B_{\frac{\rho}{2}}}$, and set
\begin{equation}\label{aaa1}
w_h(y):=v_h(y)-b_h -B_hy,
\end{equation}
for any $y\in B_1$. Proceeding  similarly as in \eqref{resc}, we rescale $F$ and $G$ around $A_h+\lambda_h B_h$,
\begin{align}\label{riscf}
& \widetilde F_h(\xi):= \frac{F( A_h +\lambda_h B_h+\lambda_h\xi)- F( A_h +\lambda_h B_h)-D F( A_h +\lambda_h B_h)\lambda_h\xi}{\lambda_h^2},\\
& \widetilde G_h(\xi):= \frac{G( A_h +\lambda_h B_h+\lambda_h\xi)- G( A_h +\lambda_h B_h)-D G( A_h +\lambda_h B_h)\lambda_h\xi}{\lambda_h^2},
\end{align}
for any $\xi\in \R^{n\times N}$. By Lemma \ref{nic}, two growth estimates on $\widetilde F_h$, $\widetilde G_h$ and their gradients hold with some constants that depend on $p,L_1,L_2,M$ (see \eqref{decay04}) and could also depend on $\rho$
through $|\lambda_h B_h |$. However, given $\rho$,  we may  choose $h_0=h_0(n,p,M,\rho)$ large enough to have 
\begin{equation}
    |\lambda_h B_h | <\frac{c(n,p,M)\lambda_h}{\rho^{\frac{n}{p}}}<1,
\end{equation} for any $h\geq h_0$. Indeed, by \eqref{24bis} the sequence $\left\{D v_h\right\}_h$ is equibounded in $L^p(B_1)$, then we have

 \begin{align} \label{Bh}|B_h|
 & \leq \frac{2^n}{\omega_n\rho^{\frac{n}{p}}}\bigg[\int_{B_{\frac{\rho}{2}}\cap\{|D v_h|\leq 1\}}|Dv_h|\,dy+\int_{B_{\frac{\rho}{2}}\cap\{|D v_h|> 1\}}|Dv_h|\,dy\bigg]\\
& \leq \frac{2^n}{\omega_n\rho^{\frac{n}{p}}}\Bigg[\bigg(\int_{B_{\frac{\rho}{2}}}|V(Dv_h)|^2\,dy\bigg)^{\frac{1}{2}}+\bigg(\int_{B_{\frac{\rho}{2}}}|V(Dv_h)|^2\,dy\bigg)^{\frac{1}{p}}\Bigg]\leq\frac{c(n,p,M)}{\rho^{\frac{n}{p}}},\end{align}

\noindent and so the constant in \eqref{R1} can be taken independently of $\rho$.

\noindent Set 
\begin{equation}
\psi_{1,h}:=\eta w_h\quad \mathrm{and}\quad \psi_{2,h}:=(1-\eta) w_h.
\end{equation}
By the uniformly strict quasiconvexity of  $\widetilde F_h$ we have 
\begin{align}\label{decay611b}
&\frac{\ell_1}{\lambda_h^2}\int_{B_s}|V(\lambda_h D w_h)|^2\,dy\notag\\
& \leq \ell_1\int_{B_t}\big(1+|\lambda_h D\psi_{1,h}|^2\big)^{\frac{p-2}{2}}|D\psi_{1,h}|^2\,dy \leq \int_{B_t}\widetilde F_h(D\psi_{1,h})\,dy\notag\\
& =\int_{B_t}\widetilde F_h( Dw_h)\,dy+\int_{B_t}\widetilde F_h( Dw_h-D\psi_{2,h})\,dy-\int_{B_t}\widetilde F_h( Dw_h)\,dy\notag\notag\\
& =\int_{B_t}\widetilde F_h( Dw_h)\,dy-\int_{B_t}\int_0^1 D\widetilde F_h( Dw_h-\theta D\psi_{2,h})D\psi_{2,h}\,d\theta\,dy.
\end{align}
We estimate separately the two addends in the right-hand side of the previous chain of inequalities. We deal with the first addend by means of a rescaling of the minimality condition of $(u,E)$.
Using the change of variable $x=x_h+ r_hy$, the fact that $G\ge 0$ and  the minimality of $(u,E)$ with respect to $(u+\varphi,E)$ for $\varphi\in W^{1,p}_0(B_{r_h}(x_h);\R^{N})$, we have
\begin{align*}
&\int_{B_1}F(Du(x_h+r_hy))dy\leq \int_{B_1}\big[F(Du(x_h+r_hy))+\mathbbm{1}_{E^*_h}G(Du(x_h+r_hy))\big]\,dy\\
& \leq \int_{B_1} \big[ F(Du(x_h+r_hy)+D\varphi(x_h+r_hy))+\mathbbm{1}_{E^*_h}G(Du(x_h+r_hy)+D\varphi(x_h+r_hy))\big]\,dy,
\end{align*}
i.e., by the definitions of $v_h$ and $w_h$, \eqref{23} and \eqref{aaa1} respectively,
\begin{align*}
& \int_{B_1}F(A_h+\lambda_h B_h+\lambda_h Dw_h)dy\\
& \leq\int_{B_1} \big[F(A_h+\lambda_hB_h+\lambda_h(Dw_h+D\psi))+\mathbbm{1}_{E^*_h}G(A_h+\lambda_hB_h+\lambda_h(Dw_h+D\psi))\,dy,
\end{align*}
for $\psi:=\frac{\varphi(x_h+r_hy)}{\lambda_h r_h}\in W^{1,p}_0(B_1;\R^{N})$.
Therefore, recalling the  definitions of $\widetilde F_h$ and $\widetilde G_h$ in \eqref{riscf}, we have that
\begin{align*}
& \int_{B_1}\widetilde F_h(Dw_h)dy
\leq\int_{B_1}\big[\widetilde F_h(Dw_h+D\psi)+\mathbbm{1}_{E^*_h}\widetilde G_h(Dw_h+D\psi)\big]\,dy\\
& +\frac{1}{\lambda_h^2}\int_{B_1}\mathbbm{1}_{E^*_h}\big[ G(A_h+\lambda_h B_h)+ DG(A_h+\lambda_h B_h)\lambda_h(D w_h+D\psi)\big]\,dy.
\end{align*}
Choosing  $\varphi \;\mathrm{such\; that}\;  \psi=-\psi_{1,h}$,  the previous inequality becomes
\begin{align}\label{decay611c}
&\int_{B_t}\widetilde  F_h(D w_h)\,dy
 \leq \int_{B_t}\big[\widetilde  F_h\big( D w_h-D\psi_{1,h}\big)
+\mathbbm{1}_{E^*_h}\widetilde G_h( D w_h-D\psi_{1,h})\big]\,dy\\
& +\frac{1}{\lambda_h^2}\int_{B_1}\mathbbm{1}_{E^*_h}\big[G(A_h+\lambda_h B_h)+ D G(A_h+\lambda_h B_h)\lambda_h(D w_h-D\psi_{1,h})\big]\,dy\notag\\
& =\int_{B_t\setminus B_s}\big[\widetilde F_h(D\psi_{2,h})+\mathbbm{1}_{E^*_h}\widetilde G_h(D\psi_{2,h})\big]\,dy\notag\\
& +\frac{1}{\lambda_h^2} \int_{B_1}\mathbbm{1}_{E^*_h}\big[G(A_h+\lambda_h B_h)+ D G(A_h+\lambda_h B_h)\lambda_h D\psi_{2,h}\big]\,dy\notag\\
&  \leq \frac{c(p,L_1,L_2,M)}{\lambda_h^2}\int_{B_t\setminus B_s}|V(\lambda_h D\psi_{2,h})|^2\,dy+c(n,p,L_2,M)\bigg[\frac{|E^*_h|}{\lambda_h^2}+\frac{1}{\lambda_h}\int_{ E^*_h} |D\psi_{2,h}|\,dy\bigg],
\end{align}
where we have used Lemma \ref{nic}, the second estimate in \eqref{(H4)}, and the fact that $|A_h+\lambda_h B_h|\leq M+1$. By applying H\"older's and Young's inequalities, we get
\begin{align*}
    \frac{1}{\lambda_h}\int_{ E^*_h} |D\psi_{2,h}|\,dy
    & \leq \frac{|E^*_h|^{\frac{p-1}{p}}}{\lambda_h^2}\bigg(\int_{E^*_h\cap(B_t\setminus B_s)}|\lambda_h D\psi_{2,h}|^p\,dy\bigg)^{\frac{1}{p}}\\
    & \leq \frac{1}{\lambda_h^2}\bigg[|E^*_h|+\int_{E^*_h\cap(B_t\setminus B_s)}|\lambda_h D\psi_{2,h}|^p\,dy\bigg]\\
    & \leq \frac{1}{\lambda_h^2}\bigg[2|E_h^*|+\int_{E^*_h\cap (B_t\setminus B_s)\cap\{|\lambda_h D\psi_{2,h}|>1\}}|\lambda D\psi_{2,h}|^p\,dy\bigg]\\
    & \leq \frac{1}{\lambda_h^2}\bigg[2|E^*_h|+\int_{B_t\setminus B_s}|V(\lambda_h D\psi_{2,h}))|^2\,dy\bigg].
\end{align*}
The previous chain of inequalities combined with \eqref{decay611c} yields
\begin{align}\label{decay611ca}
    \int_{B_1}\widetilde F_h(Dw_h)dy
\leq \frac{c(n,p,L_1,L_2,M)}{\lambda_h^2}\bigg[\int_{B_t\setminus B_s}|V(\lambda_h D\psi_{2,h})|^2\,dy+|E^*_h|\bigg].
\end{align}

\noindent Now we estimate the second addend in the right-hand side of \eqref{decay611b}. Using the upper bound on $D\widetilde F_h$ in Lemma \ref{nic}, 
\begin{align}\label{decay611d}
    &\int_{B_t}\int_0^1 D\widetilde F_h( Dw_h-\theta D\psi_{2,h})D\psi_{2,h}\,d\theta dy\\
    & \leq c(p,L_1,M)\int_{B_t\setminus B_s}\int_0^1\big(1+\lambda_h^2|Dw_h-\theta D\psi_{2,h}|^2\big)^{\frac{p-2}{2}}|Dw_h-\theta D\psi_{2,h}||D\psi_{2,h}|d\theta dy. 
\end{align}
Regarding the integrand in the latest estimate, we distinguish two cases:

\medskip

\indent\textbf{Case 1:}
$|D\psi_{2,h}|\leq|Dw_h-\theta D\psi_{2,h}|$.\\
By the definition of $V$, we have
\begin{equation*}
\big(1+\lambda_h^2|Dw_h-\theta D\psi_{2,h}|^2\big)^{\frac{p-2}{2}}|Dw_h-\theta D\psi_{2,h}||D\psi_{2,h}|\leq \lambda_h^{-2} |V(\lambda_h( Dw_h-\theta D\psi_{2,h})|^2.
\end{equation*}
\indent \textbf{Case 2:}
$|Dw_h-\theta D\psi_{2,h}|<|D\psi_{2,h}|$.\\
If $|D\psi_{2,h}|<1/\lambda_h$, using \emph{(i)} of Lemma \ref{PropV} we get
\begin{equation*}
\big(1+\lambda_h^2|Dw_h-\theta D\psi_{2,h}|^2\big)^{\frac{p-2}{2}}|Dw_h-\theta D\psi_{2,h}||D\psi_{2,h}|\leq |D\psi_{2,h}|^2\leq \lambda_h^{-2} |V(\lambda_h D\psi_{2,h})|^2.
\end{equation*}
If $|D\psi_{2,h}|\geq 1/\lambda_h$, using  again (i) of Lemma \ref{PropV} we deduce that
\begin{align}
&\big(1+\lambda_h^2|Dw_h-\theta D\psi_{2,h}|^2\big)^{\frac{p-2}{2}}|Dw_h-\theta D\psi_{2,h}||D\psi_{2,h}|\leq \\
    & \leq \lambda_h^{p-2}|Dw_h-\theta D\psi_{2,h}|^{p-1}|D\psi_{2,h}|\leq \lambda_h^{-2}|\lambda_h D\psi_{2,h}|^{p}\leq \lambda_h^{-2} |V(\lambda_h D\psi_{2,h})|^2. 
\end{align}
By combining the two previous cases, we can proceed in the estimate \eqref{decay611d} as follows:
\begin{align}\label{decay711d}
    &\int_{B_t}\int_0^1 D\widetilde F_h( Dw_h-\theta D\psi_{2,h})D\psi_{2,h}\,d\theta\,dy \\
    & \leq \frac{c(p,L_1,M)}{\lambda_h^2}\int_{B_t\setminus B_s} \big(|V(\lambda_h( Dw_h-\theta D\psi_{2,h})|^2+|V(\lambda_h D\psi_{2,h})|^2\big)\,dy\\
    & \leq \frac{c(p,L_1,M)}{\lambda_h^2}\int_{B_t\setminus B_s} \big(|V(\lambda_h Dw_h)|^2+|V(\lambda_h D\psi_{2,h})|^2\big)\,dy
\end{align}
Hence, combining \eqref{decay611b} with \eqref{decay611ca} and \eqref{decay711d}, we obtain
\begin{align*}
   &\frac{\ell_1}{\lambda_h^2}\int_{B_s}|V(\lambda_h D w_h)|^2\,dy\\ 
   &\leq \frac{c(n,p,L_1,L_2,M)}{\lambda_h^2}\bigg[\int_{B_t\setminus B_s}\big(|V(\lambda_h Dw_h)|^2+|V(\lambda_h D\psi_{2,h})|^2\big)\ dy+|E^*_h|\bigg]
\end{align*}
By the definition of $\psi_{2,h}$ and \emph{(iii)} of Lemma \ref{PropV}, we infer that
\begin{align}
&\ell_1\int_{B_s}|V(\lambda_h D w_h)|^2\,dy\\ 
   &\leq \tilde{C}\bigg[\int_{B_t\setminus B_s}\bigg(|V(\lambda_h Dw_h)|^2+\bigg|V\bigg(\lambda_h \frac{w_h}{t-s}\bigg)\bigg|^2\bigg)\, dy+|E^*_h|\bigg],
\end{align}
for some $\tilde{C}=\tilde{C}(n,p,L_1,L_2,M)$\\
 \indent By adding $\tilde{C}\int_{B_s}|V(\lambda_h D w_h)|^2\,dy$ to both sides of the previous estimate, dividing by $\ell_1+\tilde{C}$ and thanks to Lemma \ref{iterationV}, we deduce that
\begin{equation*}
   \int_{B_{\frac{\rho}{2}}} |V(\lambda_h D w_h)|^2\,dy\leq c(n,p,\ell_1,L_1,L_2,M)\bigg(\int_{B_\rho}\bigg|V\bigg(\lambda_h\frac{ w_h}{\rho}\bigg)\bigg|^2\,dy+|E^*_h|\bigg). 
\end{equation*}
Therefore, by the definition of $w_h$, we conclude that
\begin{align}\label{Caccio}
& \int_{B_{\frac{\rho}{2}}}\big|V(\lambda_h(Dv_h -(Dv_h)_{\frac{\rho}{2}})\big|^2\,dy\\
& \leq c(n,p,\ell_1,L_1,L_2,M)\Bigg[\int_{B_\rho}\bigg|V\bigg(\frac{\lambda_h\bigl(v_h -(v_h)_\rho-(Dv_h)_{\frac{\rho}{2}}\,y\bigr)}{\rho}\bigg)\bigg|^2\,dy+|E^*_h|\Bigg]\notag
\end{align}
which, by the relative isoperimetric inequality and the hypothesis of this substep, i.e. $\min\{|E^*_h|, |B_1\setminus E_h|\}=|E^*_h|$, yields the estimate \eqref{Cacciofin}.

\medskip

\noindent{\bf Substep 2.b}\,\,{\it  The case $\min\{|E^*_h|, |B_1\setminus E_h|\}=|B_1\setminus E_h|$.}
\medskip

\noindent As in the previous substep, we fix  $0<\frac{\rho}{2}<s<t<\rho<1$ and let $\eta\in C^\infty_0(B_{t})$  be a cut off function between $B_s$ and $B_t$, i.e., $0\leq \eta\leq 1$, $\eta\equiv 1$ on $B_s$ and $|\nabla \eta |\leq \frac{c}{t-s}$.
Also, we set $b_h := (v_h)_{B_\rho}$ , $B_h := (Dv_h)_{B_{\frac{\rho}{2}}}$ and define
\begin{equation}
w_h(y):=v_h(y)-b_h -B_hy, \quad\forall y\in B_1,
\end{equation}
and
\begin{equation*}
\widetilde H_h:=\widetilde F_h+
\widetilde G_h.
\end{equation*}
We remark that Lemma \ref{nic} can be applied to $\widetilde H_h$, that is
\begin{equation*}
    |\widetilde H_h(\xi)|\leq c(p,L_1,L_2,M)\big(1+|\lambda_h\xi|^2\big)^{\frac{p-2}{2}}|\xi|^2 ,\quad\forall\xi\in\R^{n\times N},
\end{equation*}
and, by the uniformly strict quasiconvexity conditions $\eqref{F2p}$ and $\eqref{G2p}$,
\begin{equation}\label{quasiH}
\int_{B_1}\widetilde H_h(\xi+D \psi)\, dx \ge \int_{B_t}\big[\widetilde{H}_h(\xi)+\tilde{\ell}\big(1+|\lambda_h D\psi|^2\big)^{\frac{p-2}{2}}|D\psi|^2\big]\,dy,\quad\forall\psi\in W^{1,p}_0(B_1;\R^{N}),
\end{equation}
 where we have set
\begin{equation*}
    \widetilde \ell= \ell_1+\ell_2.
\end{equation*}
\noindent We set again  \begin{equation}
    \psi_{1,h}:=\eta w_h\quad \mathrm{and} \quad \psi_{2,h}:=(1-\eta) w_h.
\end{equation}
By the quasiconvexity condition \eqref{quasiH} and since $\widetilde H_h(0)=0$, we have
 \begin{align}\label{decay611e} 
\frac{\widetilde{\ell}}{\lambda_h^2} \int_{B_s}&|V(\lambda_h D w_h)|^2\,dy\leq
\widetilde{\ell}\int_{B_s}\big(1+|\lambda_h D w_h|^2\big)^{\frac{p-2}{2}}|D w_h|^2\,dy\notag\\
& \leq \widetilde{\ell}\int_{B_t}\big(1+|\lambda_h D\psi_{1,h}|^2\big)^{\frac{p-2}{2}}|D\psi_{1,h}|^2\,dy \leq \int_{B_t}\widetilde H_h(D\psi_{1,h})\,dy=\int_{B_t}\widetilde H_h( Dw_h-D\psi_{2,h})\,dy\notag\\
& =\int_{B_t}\widetilde H_h( Dw_h)\,dy+\int_{B_t}\widetilde H_h( Dw_h-D\psi_{2,h})\,dy-\int_{B_t}\widetilde H_h( Dw_h)\,dy\notag\\
& =\int_{B_t}\widetilde H_h( Dw_h)\,dy-\int_{B_t}\int_0^1 D\widetilde H_h( Dw_h-\theta D\psi_{2,h})D\psi_{2,h}\,d\theta\,dy.
\end{align}
Similarly to the previous case, we estimate separately the two addends in the right-hand side of the previous chain of inequalities. Using the minimality condition \eqref{decay5aaa} for the rescaled functions $v_h$ and recalling the definition of $\tilde{H}_h$, since $Dv_h=Dw_h+B_h$, we get
\begin{align}\label{decay5aa} \int_{B_1}\widetilde H_h(Dw_h)dy
& \leq\int_{B_1}\widetilde H_h(Dw_h+D\psi)\,dy\notag\\
& + \frac{L_2}{\lambda_h^2}\int_{(B_1\setminus E_h)\cap \mathrm{supp}\psi}\big(1+|A_h+\lambda_hB_h+\lambda_hDw_h|^2\big)^{\frac{p}{2}}\,dy.
\end{align}
Choosing $\psi=-\psi_{1,h}$ as test function in \eqref{decay5aa} and using the fact that $\widetilde H_h(0)=0$, we estimate
\begin{align}\label{decay611a}
& \int_{B_t}\widetilde  H_h(D w_h)\,dy\notag\\
& \leq \int_{B_t}\widetilde  H_h( D w_h-D\psi_{1,h})\,dy+\frac{L_2}{\lambda_h^2}\int_{B_t\setminus E_h}\big(1+|A_h+\lambda_hB_h+\lambda_hDw_h|^2\big)^{\frac{p}{2}}\,dy\notag\\
& =\int_{B_t\setminus B_s}\widetilde  H_h\big( D \psi_{2,h}\big)\,dy+\frac{L_2}{\lambda_h^2}\int_{B_t\setminus E_h}\big(1+|A_h+\lambda_hB_h+\lambda_hDw_h|^2\big)^{\frac{p}{2}}\,dy\\
& \leq \frac{c(p,L_1,L_2,M)}{\lambda_h^2}\int_{B_t\setminus B_s}\big|V\big(\lambda_hD\psi_{2,h}\big)|^2\,dy+\frac{L_2}{\lambda_h^2}\int_{B_t\setminus E_h}\big(1+|A_h+\lambda_hB_h+\lambda_hDw_h|^2\big)^{\frac{p}{2}}\,dy.
\end{align}
We note that, since $|A_h+\lambda_hB_h|\leq c(M)$, for every fixed $\varepsilon>0$ there exists a constant $C=C(p,\varepsilon)$ such that
\begin{equation*}
\big(1+|A_h+\lambda_h B_h+\lambda_h Dw_h|^2\big)^{\frac p2}\leq C(p,\varepsilon) c(M)^p+(1+\varepsilon)\lambda_h^p|D w_h|^p.
\end{equation*}
Summarizing, we get
\begin{align}\label{decay611f}
&\int_{B_t}\widetilde  H_h(D w_h)\,dy\leq \frac{c(p,L_1,L_2,M)}{\lambda_h^2}\int_{B_t\setminus B_s}\big|V\big(\lambda_hD\psi_{2,h}\big)|^2\,dy\\
& +(1+\varepsilon)\frac{L_2}{\lambda_h^2}\int_{B_t}\mathbbm{1}_{_{\{|\lambda_h D w_h|\geq 1\}}}|\lambda_hD w_h|^p\,dy+c(p,L_2,M,\varepsilon)\frac{|B_1\setminus E_h|}{\lambda_h^2}\notag.
\end{align}
Now we estimate the second addend in the right-hand side of \eqref{decay611e}. Using the upper bound on $D\widetilde H_h$ in Lemma \ref{nic}, we obtain
\begin{align}\label{decay611g}
&\int_{B_t}\int_0^1 D\widetilde H_h( Dw_h-\theta D\psi_{2,h})D\psi_{2,h}\,d\theta\,dy\\
    & \leq c(p,L_1,L_2,M)\int_{B_t\setminus B_s}\int_0^1\big(1+\lambda_h^2|Dw_h-\theta D\psi_{2,h}|^2\big)^{\frac{p-2}{2}}|Dw_h-\theta D\psi_{2,h}||D\psi_{2,h}|\,dy. 
\end{align}
Proceeding exactly as in the estimate \eqref{decay711d} of the step $2.a$, we obtain
\begin{align}\label{decay711ds}
    &\int_{B_t}\int_0^1 D\widetilde H_h( Dw_h-\theta D\psi_{2,h})D\psi_{2,h}\,d\theta\,dy \\
    & \leq \frac{c(p,L_1,L_2,M)}{\lambda_h^2}\int_{B_t\setminus B_s} \big(|V(\lambda_h Dw_h)|^2+|V(\lambda_h D\psi_{2,h})|^2\big)\,dy.
\end{align}
Inserting \eqref{decay611f} and \eqref{decay711ds} in \eqref{decay611e}, we infer that
\begin{align*}
\frac{\widetilde{\ell}}{\lambda_h^2} \int_{B_s}&|V(\lambda_h D w_h)|^2\,dy\\
& \leq \frac{c(p,L_1,L_2,M)}{\lambda_h^2}\int_{B_t\setminus B_s} \big(|V(\lambda_h Dw_h)|^2+|V(\lambda_h D\psi_{2,h})|^2\big)\ dy\\
& +(1+\varepsilon)\frac{L_2}{\lambda_h^2}\int_{B_t}\mathbbm{1}_{_{\{|\lambda_h D w_h|\geq 1\}}}|\lambda_hD w_h|^p\,dy+c(p,L_2,M,\varepsilon)\frac{|B_1\setminus E_h|}{\lambda_h^2} \\
& \leq \frac{c(p,L_1,L_2,M)}{\lambda_h^2}\int_{B_t\setminus B_s} |V(\lambda_h Dw_h)|^2\,dy+\frac{c(p,M,L_1,L_2)}{\lambda_h^2}\int_{B_t\setminus B_s}\bigg|V\bigg(\lambda_h \frac{w_h}{t-s}\bigg)\bigg|^2\,dy\\
& +(1+\varepsilon)\frac{L_2}{\lambda_h^2}\int_{B_t} |V(\lambda_h Dw_h)|^2\,dy+c(p,L_2,M,\varepsilon)\frac{|B_1\setminus E_h|}{\lambda_h^2}.
\end{align*}
 Taking advantage of the  hole filling technique as in the previous case, we obtain
\begin{align*}\label{decay611bbba}
& \int_{B_t\setminus B_s} |V(\lambda_h Dw_h)|^2\ dy\\
& \leq\frac{(c(p,L_1,L_2,M) + (1+\varepsilon)L_2)}{(c(p,M,L_1,L_2)+\widetilde\ell)}\int_{B_t} |V(\lambda_h Dw_h)|^2\ dy\\
& +c(p,M,L_1,L_2)\int_{B_t\setminus B_s}\bigg|V\bigg(\lambda_h \frac{w_h}{t-s}\bigg)\bigg|^2\,dy+c(p,L_2,M,\varepsilon)\frac{|B_1\setminus E_h|}{\lambda_h^2}.
\end{align*}
 The assumption \eqref{H} implies that there exists $\varepsilon=\varepsilon(p,\ell_1,\ell_2,L_2)>0$ such that $\frac{(1+\varepsilon)L_2}{\ell_1+\ell_2}<1$. Therefore
  we have 
  \begin{equation}
      \frac{c + (1+\varepsilon)L_2}{c+\widetilde\ell}= \frac{c + (1+\varepsilon) L_2}{c+\ell_1+\ell_2}<1.
  \end{equation}
So, by virtue of Lemma \ref{iterationV}, from the previous estimate we deduce that
\begin{equation*}
\int_{B_{\frac{\rho}{2}}} |V(\lambda_h D w_h)|^2\,dy\leq c(n,p,\ell_1,\ell_2,L_1L_2,M)\bigg(\int_{B_\rho}\bigg|V\bigg(\lambda_h\frac{w_h}{\rho}\bigg)\bigg|^2\,dy+|B_1\setminus E_h|\bigg).
\end{equation*}
By definition of $w_h$ and the relative isoperimetric inequality, since  $|B_1\setminus E_h|=\min\{|E^*_h|, |B_1\setminus E_h|\}$, we get the estimate \eqref{Cacciofin}.
\medskip

\textbf{Step 3.}
\emph{$v$ solves a linear system in $B_1$}.
 
Let us divide the proof into two cases, depending on which one is the smallest between $|E^*_h|$ and $|B_1\setminus E_h|$.

\medskip
We divide the proof in two substeps.

\medskip

\noindent{\bf Substep 3.a}\,\,{\it  The case $\min\{|E^*_h|, |B_1\setminus E_h|\}=|E^*_h|$.}
\medskip
\noindent
We claim that $v$ solves the linear system
 \begin{equation*}
     \int_{B_{1}}D^2
F(A) DvD\psi \,dy=0,
 \end{equation*}
 for all $\psi \in C^1_0(B_1;\R^N)$.
Since $v_h$ satisfies  \eqref{29}, we have that
\begin{equation}\label{elext}
0 \leq {\cal
I}_h(v_h+s\psi)-{\cal I}_h(v_h)+\frac{1}{\lambda_h}\int_{B_1}\mathbbm{1}_{E^*_h} D G(A_h)sD\psi\,dy ,
\end{equation}
for every $\psi \in C^{1}_{0}(B_{1};\R^N)$ and  $s\in (0,1)$. 
Dividing by $s$ and passing to the limit as $s\to 0$, by  the definition of ${\cal I}_{h}$, we get  (see \cite{CEL} or \cite{CFP1})
\begin{align}\label{47}
0
&\leq\frac{1}{\lambda_h}\int_{B_{1}}( DF( A_h+\lambda_hDv_h) - DF(A_h))D\psi\, dy\notag\\
& +\frac{1}{\lambda_h}\int_{B_{1}}\mathbbm{1}_{E^*_h} DG( A_h+\lambda_hDv_h) D\psi \,dy.
\end{align}
We partition the unit ball as follows:
\begin{equation}
B_{1}=\mathbf{B}_h^+\cup \mathbf{B}_h^-=\{y\in B_{1}:\, \lambda_h|Dv_h|>1\}\cup\{y\in B_{1}:\, \lambda_h|Dv_h|\leq 1\} .
\end{equation}
By \eqref{24bis}, we get

\begin{equation}\label{48g2}
|\mathbf{B}_h^+|\leq \int_{\mathbf{B}_h^+}\lambda_h^{p}|Dv_h|^{p}\, dy\leq \lambda_h^{p}\int_{B_1}|Dv_h|^{p}\, dy\leq c(n,p,M)\lambda_h^{p} .
\end{equation}
We rewrite \eqref{47} as follows:
    \begin{align}
    \label{L1}
0
&\leq \frac{1}{\lambda_h}\int_{\mathbf{B}_h^+}( DF( A_h+\lambda_hDv_h) - DF( A_h))D\psi\, dy\notag\\
& +\int_{\mathbf{B}_h^-}\int_0^1\left( D^2F( A_h+t\lambda_hDv_h) - D^2F( A)\right)dt Dv_h D\psi\, dy\notag\\
& +\int_{\mathbf{B}_h^-}D^2 F(A)Dv_hD\psi\,dy+\frac{1}{\lambda_h}\int_{B_{1}}\mathbbm{1}_{E^*_h} DG( A_h+\lambda_hDv_h) D\psi \,dy.
\end{align}

By growth condition in  \eqref{(H4)} and
 H\"older's inequality, we get
\begin{align}\label{49}
&\frac{1}{\lambda_h}\left|\int_{\mathbf{B}_h^+}(D
F(A_h+\lambda_hDv_h)-D
F(A_h))D\psi \,dy\right|\\
&\leq c(p,L_1,M,D\psi)\bigg[\frac{|\mathbf{B}_h^+|}{\lambda_h}+\lambda_h^{p-2}\int_{\mathbf{B}_h^+}|Dv_h|^{p-1}\, dy\bigg]\\
&\leq c(n,p,L_1,M,D\psi)\Bigg[\lambda_h^{p-1}+\lambda_h^{p-1}\left(\int_{\mathbf{B}_1}|Dv_h|^p\, dy\right)^{\frac{p-1}{p}}\bigg(\frac{|\mathbf{B}_h^+|}
{\lambda_h^{p}}\bigg)^{\frac{1}{p}}\Bigg]\leq c(n,p,L_1,M,D\psi)\lambda_h^{p-1},
\end{align}
thanks to  \eqref{decay04}, \eqref{24bis} and \eqref{48g2}. Thus
\begin{equation}\label{50}
\lim_{h\to \infty}\frac{1}{\lambda_h}\bigg|\int_{\mathbf{B}_h^+}(D
F(A_h+\lambda_hDv_h)-D
F(A_h))D\psi \,dy\bigg|=0.
\end{equation}
By \eqref{decay04} and the definition of $\mathbf{B}_h^-$ we have that $|A_h+\lambda_hDv_h|\leq M+1$ on $\mathbf{B}_h^-$. Hence we estimate 
\begin{align}\label{448}
&\bigg|\int_{\mathbf{B}_h^-}\int_0^1\left (D^2\notag
F(A_h+t\lambda_hDv_h)-D^2F(A)\right)dt Dv_hD\psi \,dy\bigg|\\
& \leq\int_{\mathbf{B}_h^-}\left|\int_0^1\left(D^2
F(A_h+t\lambda_hDv_h)-D^2F(A)\right)\,dt\right| |Dv_h||D\psi| \,dy\notag\\
& \leq \left(\int_{\mathbf{B}_h^-}\left|\int_0^1\left(D^2
F(A_h+t\lambda_hDv_h)-D^2F(A)\right)\,dt\right|^{\frac{p}{p-1}}dy\right)^{\frac{p-1}{p}} \norm{Dv_h}_{L^{p}(B_1)}\norm{D\psi}_{L^\infty(B_1)}\notag\\
& \leq c(n,p,M,D\psi)\left(\int_{\mathbf{B}_h^-}\left|\int_0^1\left(D^2
F(A_h+t\lambda_hDv_h)-D^2F(A)\right)\,dt\right|^{
\frac{p}{p-1}}dy\right)^{\frac{p-1}{p}},
\end{align}
where we have used \eqref{24bis}. Since, by \eqref{25},
$\lambda_hDv_h\to 0$ a.e. in $B_{1}$, the uniform continuity of  $D^2F$ on bounded sets implies that
\begin{equation}\label{L2}
\lim_{h\rightarrow\infty}\bigg|\int_{\mathbf{B}_h^-}\int_0^1\left (D^2
F(A_h+t\lambda_hDv_h)-D^2F(A)\right)dt Dv_hD\psi \,dy\bigg|=0.
\end{equation}
 Note that \eqref{48g2} yields that $\mathbbm{1}_{\mathbf{B}_h^-}\to \mathbbm{1}_{B_1}$ in $L^r(B_1)$, for every $r<\infty$. Therefore, by the weak convergence of $D v_h$ to $D v$ in $L^p(B_1)$, it follows that
  \begin{equation}\label{L3}
\lim_{h\rightarrow\infty}\int_{\mathbf{B}_h^-}D^2F(A) Dv_hD\psi\,dy=\int_{B_1}D^2F(A) DvD\psi\,dy.
\end{equation}
By growth condition  \eqref{(H4)}, we deduce 
\begin{align*}
& \frac{1}{\lambda_h}\left|\int_{B_1}\mathbbm{1}_{E^*_h}[D_\xi
G(A_h+\lambda_hDv_h)D\psi \,dy\right|
\leq \frac{c(p,L_2)}{\lambda_h}\int_{B_1}\mathbbm{1}_{E^*_h}\big( 1+|A_h+\lambda_hDv_h|^2\big)^{\frac{p-1}{2}}|D\psi| \,dy\\
&\leq c(p,L_2,M,D\psi)\bigg[\frac{1}{\lambda_h}|E^*_h|+\lambda_h^{p-2}\int_{ E^*_h}|D v_h|^{p-1}\, dy\bigg]\\
&\leq c(p,L_2,M,D\psi)\bigg[\frac{1}{\lambda_h}|E^*_h|+\lambda_h^{p-2+\frac{2}{p}}\left(\int_{B_1}|Dv_h|^p\, dy\right)^{\frac{p-1}{p}}\bigg(\frac{|E^*_h|}
{\lambda_h^2}\bigg)^{\frac{1}{p}}\bigg]\\
&\leq  c(n,p,L_2,M,D\psi)\bigg[\frac{1}{\lambda_h}|E^*_h|+\lambda_h^{p-2-\frac{2}{p}}\bigg(\frac{|E^*_h|}
{\lambda_h^2}\bigg)^{\frac{1}{p}}\bigg],
\end{align*}
where we have used \eqref{decay04}  and  \eqref{24bis}. Since $\min\{|E^*_h|, |B_1\setminus E_h|\}=|E^*_h|$,  by   \eqref{25ter}, we have
\begin{equation*}
    \lim_{h\rightarrow\infty} \frac{|E^*_h|
}{\lambda_h^2}=0,
\end{equation*}
and so
\begin{equation}\label{50bis}
\lim_{h\to \infty}\frac{1}{\lambda_h}\int_{B_1}\mathbbm{1}_{E^*_h}D
G(A_h+\lambda_hDv_h)D\psi \,dy=0 .
\end{equation}
By \eqref{50}, \eqref{L2}, \eqref{L3} and \eqref{50bis},
passing to the limit as $h\to \infty$ in \eqref{L1}, we get
\begin{equation}
\int_{B_{1}}D
F(A)  DvD\psi \,dy\geq 0.
\end{equation}
Furthermore, plugging $-\psi$ in place of $\psi$, we  get
\begin{equation*}\label{prob000} \int_{B_{1}}D
F(A) DvD\psi \,dy=0,\end{equation*}
i.e. $v$ solves a linear system with constant coefficients.

\medskip

\noindent{\bf Substep 3.b}\,\,{\it  The case $\min\{|E^*_h|, |B_1\setminus E_h|\}=|B_1\setminus E_h|$.}

\medskip

\noindent We claim that $v$ solves the linear system
\begin{equation*}
   \int_{B_{1}}D^2
(F+G)(A) DvD\psi \,dy=0, 
\end{equation*}
for all $\psi\in C^1_0(B_1;\R^N)$.
Dividing by $s$ and passing to the limit as $s\to 0$, by  the definition of ${\cal H}_{h}$ we get  (see \cite{CEL} or \cite{CFP1})
\begin{align}
\label{eqeq1}
& 0 \leq\frac{1}{\lambda_h}\int_{B_1}\big[D (F+G)(A_h+\lambda_hDv_h)D\psi- D(F+G)(A_h)D\psi\big] dy\\
& +c(p,L_2,M)\bigg[\frac{1}{\lambda_h}\int_{B_1\setminus E_h}|D\psi|dy+
\int_{B_1\setminus E_h}\lambda_h^{p-2}|Dv_h|^{p-1}|D\psi|\,dy\bigg].
\end{align}
As before, we partition  $B_1$ as follows:
\begin{equation}
B_{1}=\mathbf{B}_h^+\cup \mathbf{B}_h^-=\{y\in B_{1}:\,\lambda_h|Dv_h|>1\}\cup\{y\in B_{1}:\,\lambda_h|Dv_h|\leq 1\}. 
\end{equation}
We rewrite \eqref{eqeq1} as
\begin{align}
    \label{L4}
0
& \leq \frac{1}{\lambda_h}\int_{\mathbf{B}_h^+}( D(F+G)( A_h+\lambda_hDv_h) - D(F+G)( A_h))D\psi\, dy\\
& +\frac{1}{\lambda_h}\int_{\mathbf{B}_h^-}( D(F+G)( A_h+\lambda_hDv_h) - D(F+G)( A_h))D\psi\, dy\notag\\
& +c(p,L_2,M)\bigg[\frac{1}{\lambda_h}\int_{B_1\setminus E_h}|D\psi|dy+
\int_{B_1\setminus E_h}\lambda_h^{p-2}|Dv_h|^{p-1}|D\psi|\,dy\bigg]\notag.
\end{align}
Arguing as in \eqref{50}, we obtain that
\begin{align}\label{49a}
\lim_{h\to \infty}\frac{1}{\lambda_h}\left|\int_{\mathbf{B}_h^+}(D
(F+G)(A_h+\lambda_hDv_h)-D
(F+G)(A_h))D\psi \,dy\right|=0  ,
\end{align}
and, as in \eqref{L2} and \eqref{L3},
\begin{align}\label{52a}
& \lim_{h\rightarrow\infty}\frac{1}{\lambda_h}\int_{\mathbf{B}_h^-}[D
(F+G)( A_h+\lambda_hDv_h)-D
(F+G)(A_h)]D\psi \,dy\notag\\
& =\int_{B_{1}}D
(F+G)( A) DvD\psi \,dy.
\end{align}
\bigskip
Moreover, we have that
\begin{align*}
& \frac{1}{\lambda_h}\int_{B_1\setminus E_h}|D\psi|dy+
\int_{B_1\setminus E_h}\lambda_h^{p-2}|Dv_h|^{p-1}|D\psi|\,dy\\
& \leq c(p,D\psi)\Bigg[\frac{|B_1\setminus E_h|}{\lambda_h}
+ \lambda_h^{p-2+\frac{2}{p}}\left(\int_{\mathbf{B}_1}|Dv_h|^p\, dy\right)^{\frac{p-1}{p}}\bigg(\frac{|B_1\setminus E_h|}
{\lambda_h^2}\bigg)^{\frac{1}{p}}\Bigg]\\
& \leq c(n,p,D\psi)\Bigg[\frac{|B_1\setminus E_h|}{\lambda_h}
+ \lambda_h^{p-2+\frac{2}{p}}\bigg(\frac{|B_1\setminus E_h|}
{\lambda_h^2}\bigg)^{\frac{1}{p}}\Bigg],
\end{align*}
where we used \eqref{24bis}.
Since  $\min\{|E^*_h|, |B_1\setminus E_h|\}=|B_1\setminus E_h|$,  by  \eqref{25ter}, we have
\begin{equation*}
\lim_{h\rightarrow\infty} \frac{|B_1\setminus E_h|
}{\lambda_h^2}=0,
\end{equation*}
and we obtain
\begin{equation}\label{52b}
\lim_{h\rightarrow\infty} \left[\frac{1}{\lambda_h}\int_{B_1\setminus E_h}|D\psi|dy+
\int_{B_1\setminus E_h}\lambda_h^{p-2}|Dv_h|^{p-1}|D\psi|\,dy\right]=0.
\end{equation}
By  \eqref{49a}, \eqref{52a} and \eqref{52b}, passing to the limit as $h\to\infty$ in \eqref{L4} we conclude that
\begin{equation}
 \int_{B_{1}}D^2
(F+G)(A)  DvD\psi \,dy\geq 0   
\end{equation}
and, with $-\psi$ in place of $\psi$, we finally get
\begin{equation*}
\int_{B_{1}}D^2
(F+G)(A) DvD\psi \,dy=0,\end{equation*}
as claimed.

\bigskip

By Proposition \ref{regu} and the theory of linear
systems (see \cite[ Theorem 2.1 and Chapter 3]{Gia}), we deduce in both cases that
 $v \in C^\infty$ and  there exists a constant $\tilde{c}=\tilde{c}(n,N,p,\ell_1,\ell_2,L_1,L_2)>0$ such that 
\begin{equation}\label{stimafofu}\medint_{B_{\tau}}| D v-(Dv)_{\tau}|^2\leq \tilde{c} \tau^2 \medint_{B_{\frac 12}}| D v-(Dv)_{\frac 12}|^2\,dx,
\end{equation}
for any $\tau\in\big(0,\frac{1}{2}\big)$. Moreover, by Proposition \ref{regu} again,
\begin{equation}
    \medint_{B_{\frac 12}}| D v-(Dv)_{\frac 12}|^2\,dx\leq \sup_{B_{\frac 12}}|Dv|^2\leq \tilde{c}\bigg(\medint_{B_{1}}| D v|^p\,dx\bigg)^{2/p}.
\end{equation}
\noindent 
Observing that
\begin{equation}
\norm{Dv}_{L^p(B_1)}\leq \limsup_h\norm{Dv_h}_{L^p(B_1)}\leq c(n,p),
\end{equation}
it  follows that
\begin{equation}
\label{eqq2}
\medint_{B_{\tau}}| D v-(Dv)_{\tau}|^2\leq \overline{C}\tau^2,
\end{equation}
for some fixed $\overline{C}=\overline{C}(n,N,p,\ell_1,\ell_2,L_1,L_2)$.

\medskip

\noindent \textbf{Step 4.} \emph{An estimate for the perimeters.}\\
\medskip
\indent Our aim is to show that there exists a constant $c=c(n,p,L_2,\Lambda,M)>0$ such that
\begin{align}\label{StimaPer}
P(E_h,B_\tau)
\leq c\bigg[\frac{1}{\tau} P(E_h,B_1)^{\frac{n}{n-1}}+r_h \tau^n+r_h\lambda_h^{p}\bigg].
\end{align}

\noindent By the minimality of $(u, E)$ with respect to $(u,\widetilde E)$, where $\widetilde E$ is a set of finite perimeter such that $\widetilde E\Delta E\Subset B_{r_h}(x_h)$ and $B_{r_h}(x_h)$ are the balls of the contradiction argument,  we get
\begin{equation*}
    \int_{B_{r_h}(x_h)} \mathbbm{1}_E G(Du)+
\Phig(E;B_{r_h}(x_h))
\leq \int_{B_{r_h}(x_h)} \mathbbm{1}_{\widetilde E} G(Du)+ \Phig(\widetilde E;B_{r_h}(x_h)).
\end{equation*}
Using the change of variable $x=x_h+r_h y$ and dividing by $r_h^{n-1}$, we have
\begin{equation}\label{min2}
r_h\int_{B_1}\mathbbm{1}_{E_h} G(A_h+\lambda_hDv_h)dy+ \Phig_h(E_h;B_1)\leq r_h\int_{B_1}\mathbbm{1}_{{\widetilde E}_h} G(A_h+\lambda_hDv_h)dy+ \Phig_h(\widetilde E_h;B_1),
\end{equation}
where
\begin{equation*}
    \Phig_h(E_h;V):=\int_{V\cap \partial^*E_h}\Phi(x_h+r_hy,\nu_{E_h}(y))\,d\mathcal{H}^{n-1}(y),
\end{equation*}
for every Borel set $V\subset\Omega$. Assume first that $\min\{ |B_1\setminus E_h|, |E^*_h|\}=|B_1\setminus E_h|$.
Choosing $\widetilde E_h:=E_h\cup B_\rho$, we get
\begin{equation}
\label{eqq3}
    \Phig_h{(E_h;B_1)}\leq r_h\int_{B_1}\mathbbm{1}_{B_\rho} G(A_h+\lambda_hDv_h)dy+ \Phig_h(\widetilde E_h;B_1).
\end{equation}
By the coarea formula, the relative isoperimetric inequality, the choice of the representative $E^{(1)}_h$ of $E_h$, which is a Borel set, we get
\begin{equation*}
\int_\tau^{2\tau}\mathcal{H}^{n-1}(\partial B_\rho\setminus E_h)\,d\rho\leq|B_1\setminus E_h|\leq c(n) P(E_h,B_1)^{\frac{n}{n-1}}.
\end{equation*}
Therefore, thanks to Chebyshev's inequality, we may choose $\rho\in (\tau,2\tau)$, independent of $n$, such that, up to subsequences, it holds 
\begin{equation}
\label{eqq8}
   \mathcal{H}^{n-1}(\partial^*E_h\cap \partial B_{\rho})=0\quad\text{and}\quad\mathcal{H}^{n-1}(\partial B_\rho\setminus E_h)\leq \frac{c(n)}{\tau} P(E_h,B_1)^{\frac{n}{n-1}}.
\end{equation}
We remark that Proposition \ref{PerimetroAnisotropoUnione} holds also for $\Phig_h$. Thus, thanks to the choice of $\rho$, being $\mathcal{H}^{n-1}(\partial^*E_h\cap \partial B_{\rho})=0$,  we have that 
 \begin{align*}
     \Phig_h(\widetilde E_h;B_1)
    & =\Phig_h(E_h;B_\rho^{(0)})+\Phig_h(B_{\rho};E_h^{(0)})+\Phig_h(E_h;\{\nu_{E_h}=\nu_{B_\rho}\})\\
    & =\Phig_h(E_h;B_1\setminus \overline{B_\rho})+\Phig_h(B_{\rho};E_h^{(0)}).
 \end{align*}
 By the choice of the representative of $E_h$ (see Remark \ref{MTB}), taking into account \eqref{nondege} and the inequality in \eqref{eqq8}, it follows that
\begin{align}
\label{eqq6}
    \Phig_h(\widetilde E_h;B_1)
    & \leq \Phig_h(E_h;B_1\setminus \overline{B_\rho})+\Lambda\mathcal{H}^{n-1}(\partial B_{\rho}\cap E_h^{(0)})\\
    & \leq\Phig_h(E_h;B_1\setminus \overline{B_\rho})+\Lambda\mathcal{H}^{n-1}(\partial B_{\rho}\setminus E_h).\\
    & \leq \Phig_h(E_h;B_1\setminus \overline{B_\rho})+\Lambda\frac{c(n)}{\tau} P(E_h,B_1)^{\frac{n}{n-1}}.
\end{align}
On the other hand, by \eqref{nondege} and the additivity of the measure $\Phig_h(E_h,\cdot)$ it holds that
\begin{equation}
\label{eqq5}
    \frac{1}{\Lambda}P(E_h,B_\tau)\leq \Phig_h(E_h;B_\tau)\leq \Phig_h(E_h;B_1)-\Phig_h(E_h;B_1\setminus \overline{B}_\rho),
\end{equation}
since $\rho>\tau$.
Combining \eqref{eqq3}, \eqref{eqq6} and \eqref{eqq5}, we obtain
\begin{align}\label{esper}
\frac{1}{\Lambda}P(E_h,B_\tau)
& \leq \Phig_h(E_h;B_1)-\Phig_h(E_h;B_1\setminus \overline{B}_\rho)\\
&  \leq \Phig_h(\widetilde E_h;B_1)+r_h\int_{B_1}\mathbbm{1}_{B_\rho} G(A_h+\lambda_hDv_h)dy-\Phig_h(E_h;B_1\setminus \overline{B}_\rho)\\
& \leq \Lambda\frac{c(n)}{\tau} P(E_h,B_1)^{\frac{n}{n-1}}+r_h\int_{B_1}\mathbbm{1}_{B_\rho} G(A_h+\lambda_hDv_h)dy\\
& \leq \Lambda\frac{c(n)}{\tau} P(E_h,B_1)^{\frac{n}{n-1}}+c(p,L_2)r_h\int_{B_{2\tau}}(1+|A_h+\lambda_hDv_h|^2)^{\frac{p}{2}}\,dy\notag\\
& \leq \Lambda\frac{c(n)}{\tau} P(E_h,B_1)^{\frac{n}{n-1}}+c(n,p,L_2,M)r_h \tau^n+c(p,L_2)r_h\lambda_h^{p}\int_{B_{2\tau}}|D v_h|^p\,dy\notag\\
& \leq \Lambda\frac{c(n)}{\tau} P(E_h,B_1)^{\frac{n}{n-1}}+c(n,p,L_2,M)r_h \tau^n+c(n,p,L_2)r_h\lambda_h^{p},
\end{align}
where we used \eqref{24bis}. The previous estimate leads to \eqref{StimaPer}. We reach the same conclusion if $\min\{ |B_1\setminus E_h|, |E^*_h|\}=| E^*_h|$, choosing $\widetilde E_h=E_h\setminus B_\rho$ as a competitor set.

\medskip
\noindent \textbf{Step 5.} \emph{Higher integrability of $v_h$.}\\
\medskip
We need to prove that there exist two positive constants $C$ and $\delta$ depending on $n,p,\ell_1,\ell_2,L_1,L_2$ such that for every $B_r\subset B_1$ it holds
\begin{align}
\label{MaggSomm}
    \medint_{B_{\frac{r}{2}}}|V(\lambda_h D v_h)|^{2(1+\delta)}\,dy\leq C\Bigg[\bigg(\medint_{B_1}|V(\lambda_h Dv)|^2\,dy\bigg)^{1+\delta}+1\Bigg].
\end{align}
We remark that 
\begin{equation*}
    |F_h(\xi)|+|G_h(\xi)|\leq \frac{c(p,L_1,L_2,M)}{\lambda_h^2}|V(\lambda_h\xi)|^2,\quad\forall \xi\in\R^{n\times N},
\end{equation*}
and
\begin{equation*}
    \int_{B_1}F_h(D\phi)\,dy\geq \frac{\ell_1}{\lambda_h^2}\int_{B_1}|V(\lambda_h D\phi)|^2\,dy,\quad\forall\phi\in C^1_c(B_1,\R^N)
\end{equation*}
Let $r>0$ be such that $B_{3r}\subset B_1$, $\frac{r}{2}<s<t<r$ and $\eta\in C^1_c(B_t)$ be such that $0\leq\eta\leq 1$, $\eta=1$ on $B_s$, $|D\eta|\leq \frac{c}{t-s}$, for some positive constant $c$. We define
\begin{equation}
    \phi_1:=[v_h-(v_h)_r]\eta,\quad \phi_2:=[v_h-(v_h)_r](1-\eta).
\end{equation}
We deal with the case $\min\{|E^*_h|, |B_1\setminus E_h|\}=|E^*_h|$, the other one is similar. Using the minimality relation \eqref{29} and the usual growth conditions, we get
\begin{align*}
    & \frac{\ell_1}{\lambda_h^2}\int_{B_t}|V(\lambda_h D\phi_1)|^2\,dy
    \leq \int_{B_t}F_h(D\phi_1)\,dy\\
    & =\int_{B_t}F_h(D v_h)\,dy+\int_{B_t\setminus B_s}[F_h(D v_h-D\phi_2)-F_h(D v_h)]\,dy\\
    & \leq \mathcal{I}_h(v_h)+\int_{B_t\setminus B_s}[F_h(D v_h-D\phi_2)-F_h(D v_h)]\,dy\\
    & \leq \mathcal{I}_h(\phi_2)+\int_{B_t\setminus B_s}[F_h(D v_h-D\phi_2)-F_h(D v_h)]\,dy+\frac{1}{\lambda_h^2}\int_{B_t\cap E^*_h}D G(A_h)|D\phi_1|\,dy\\
    & \leq \frac{c(p,L_1,L_2,M)}{\lambda_h^2}\bigg[\int_{B_t\setminus B_s}\big[|V(\lambda_h D\phi_2)|^2+|V(\lambda_h D\phi_1)|^2+|V(\lambda_h Dv_h)|^2\big]\,dy\\
    & +\frac{1}{\lambda_h}\int_{B_t\cap E^*_h}|D\phi_1|\,dy\bigg].
\end{align*}
\medskip
By the properties of $V$, it holds that
\begin{equation*}
    |\xi|\leq C(p)\big(1+|V(\xi)|^{\frac{2}{p}}\big),\quad\forall\xi\in\R^{n\times N}.
\end{equation*}
Thus it follows
\begin{align*}
    \frac{1}{\lambda_h^2}\int_{B_t\cap E^*_h}|\lambda_h D\phi_1|\,dy & \leq \frac{c(p)}{\lambda_h^2}\bigg[|E^*_h\cap B_t|+\int_{B_t\cap E_h^*}V(|\lambda_h D\phi_1|)^{\frac{2}{p}}\,dy\bigg]\\
    & \leq \frac{c(p)}{\lambda_h^2}\bigg[c(\varepsilon)|E^*_h\cap B_t|+\varepsilon\int_{B_t\cap E_h^*}|V(\lambda_h D\phi_1)|^{2}\,dy\bigg].
\end{align*}
Combining the previous two chains of inequalities, we deduce that
\begin{align*}
    & \frac{\ell_1}{\lambda_h^2}\int_{B_t}|V(\lambda_h D\phi_1)|^2\,dy\\
    & \leq \frac{c(p,L_1,L_2,M)}{\lambda_h^2}\bigg[\int_{B_t\setminus B_s}\big[|V(\lambda_h D\phi_2)|^2+|V(\lambda_h D\phi_1)|^2+|V(\lambda_h Dv_h)|^2\big]\,dy\\
    & +c(\varepsilon)|E^*_h\cap B_t|+\varepsilon\int_{B_t\cap E_h^*}|V(\lambda_h D\phi_1)|^{2}\,dy\bigg].
    \end{align*}
    Choosing $\varepsilon$ sufficiently small, we absorb the last integral to the left-hand side 
    \begin{align}
    \label{eqqq1}
    & \frac{1}{\lambda_h^2}\int_{B_t}|V(\lambda_h D\phi_1)|^2\,dy\\
    & 
    \leq \frac{c(p,\ell_1,L_1,L_2,M)}{\lambda_h^2}\bigg[\int_{B_t\setminus B_s}\big[|V(\lambda_h D\phi_2)|^2+|V(\lambda_h D\phi_1)|^2+|V(\lambda_h Dv_h)|^2\big]\,dy+|E^*_h\cap B_t|\bigg].
    \end{align}
    By \emph{(ii)} of Lemma \ref{PropV}, it follows 
\begin{align}
    & \int_{B_s}|V(\lambda_h D v_h)|^2\,dy\\
    & 
    \leq c(p,\ell_1,L_1,L_2,M)\Bigg[\int_{B_t\setminus B_s}|V(\lambda_h D v_h)|^2\,dy+\int_{B_t\setminus B_s}\bigg|V\bigg(\lambda_h\frac{v_h-(v_h)_r}{t-s}\bigg)\bigg|^2\,dy+|E^*_h\cap B_t|\Bigg].
    \end{align}
    By applying the hole-filling technique, we add 
$c(p,\ell_1,L_1,L_2,M)\int_{B_s}|V(\lambda_h D v_h)|^2\,dy$, and we get
    \begin{align}
    & \int_{B_s}|V(\lambda_h D v_h)|^2\,dy\\
    & 
    \leq \frac{c(p,\ell_1,L_1,L_2,M)}{c(p,\ell_1,L_1,L_2,M)+1}\Bigg[\int_{B_t}|V(\lambda_h D v_h)|^2\,dy+\int_{B_t\setminus B_s}\bigg|V\bigg(\lambda_h\frac{v_h-(v_h)_r}{t-s}\bigg)\bigg|^2\,dy+|E^*_h\cap B_t| \Bigg ].
    \end{align}
    Now we can apply Lemma \ref{iterationV} and derive
    \begin{align}
    \int_{B_{r/2}}|V(\lambda_h D v_h)|^2\,dy
    \leq c(p,\ell_1,L_1,L_2,M)\Bigg[\int_{B_r}\bigg|V\bigg(\lambda_h\frac{v_h-(v_h)_r}{r}\bigg)\bigg|^2\,dy+|E^*_h\cap B_r|\Bigg].
    \end{align}
    Finally, by H\"older's inequality and Theorem \ref{SoPoV} we gain
    \begin{align*}
    \medint_{B_{r/2}}|V(\lambda_h D v_h)|^2\,dy
    & \leq c(p,\ell_1,L_1,L_2,M)\Bigg\{\Bigg[\int_{B_r}\bigg|V\bigg(\lambda_h\frac{v_h-(v_h)_r}{r}\bigg)\bigg|^{2(1+\sigma)}\,dy\Bigg]^{\frac{1}{1+\sigma}}+|B_r|\Bigg\}\\
    & \leq c(p,\ell_1,L_1,L_2,M)\bigg\{\Bigg[\medint_{B_{3r}}|V(\lambda_h D v_h)|^{\alpha}\,dy\bigg]^{\frac{1}{2\alpha}}+|B_r|\Bigg\}.
    \end{align*}
    We conclude the proof  by applying Theorem 6.6 in \cite{gi}.

\medskip
\noindent \textbf{Step 6.} \emph{ Conclusion.}
By the change of variable $x=x_h+r_h y$, \emph{(v)} of Lemma \ref{PropV} and the Caccioppoli inequality in \eqref{Cacciofin}, for every $0<\tau<\frac{1}{4}$ we have
\begin{align*}
&\limsup_{h\to \infty}\frac{U_*(x_h,\tau r_h)}{\lambda^2_h}\\
& \leq\limsup_{h\to \infty} \medint_{B_{\tau r_h}(x_0)}\bigl|V(Du)- V\bigl((Du)_{x_0,r}\bigr)\bigr|^2\,dx+\limsup_{h\to \infty} \frac{P(E, B_{\tau r_h}(x_h))}{\lambda_h^2\tau^{n-1}r_h^{n-1}}+ \limsup_{h\to \infty} \frac{\tau r_h}{\lambda_h^2}\\
&\leq \limsup_{h\to \infty}\frac{1}{\lambda_h^2}\medint_{B_{\tau}}\big|V(\lambda_h Dv_h+A_h)-V\big(A_h+\lambda_h(Dv_h)_{\tau }\big)\big|^2\,dy+ \limsup_{h\to \infty}\frac{P(E_h, B_\tau )}{\lambda_h^2\tau^{n-1}}+\tau \\
& \leq \limsup_{h\to \infty}\frac{c(n,p)}{\lambda_h^2}\medint_{B_{\tau}}\big|V(\lambda_h \big(Dv_h-(Dv_h)_{\tau }\big)\big|^2\,dy+ \limsup_{h\to \infty}\frac{P(E_h, B_\tau )}{\lambda_h^2\tau^{n-1}}+\tau\\
& \leq c(n,p,\ell_1,\ell_2,L_1,L_2,\Lambda,M) \Bigg\{\limsup_{h\to \infty}\frac{1}{\lambda_h^2}\medint_{B_{2\tau}}\bigg|V\bigg(\frac{\lambda_h\bigl(v_h -(v_h)_{2\tau}-(Dv_h)_{\tau}\,y\bigr)}{2\tau}\bigg)\bigg|^2\,dy\\
& +\frac{1}{\tau^n}\limsup_{h\to \infty} \frac{P(E_h,B_1)^{\frac{n}{n-1}}}{\lambda_h^2}+ \frac{1}{\tau^{n-1}}\limsup_{h\to \infty}\left(\frac{r_h \tau^n}{\lambda_h^2}+\frac{r_h}{\lambda_h^2}\lambda_h^p\right)+\tau\Bigg\},
\end{align*}
where we have used \eqref{24bis} and estimate \eqref{esper}.\\
Now we want to prove that
\begin{align*}
&\limsup_{h\to \infty}\frac{1}{\lambda_h^2}\int_{B_{2\tau}}\bigg|V\bigg(\frac{\lambda_h\bigl(v_h -(v_h)_{2\tau}-(Dv_h)_{\tau}\,y\bigr)}{2\tau}\bigg)\bigg|^2\,dy\\
& =\limsup_{h\to \infty}\frac{1}{\lambda_h^2}\int_{B_{2\tau}}\bigg|V\bigg(\frac{\lambda_h\bigl(v -(v)_{2\tau}-(Dv)_{\tau}\,y\bigr)}{2\tau}\bigg)\bigg|^2\,dy\leq \int_{B_{2\tau}}\frac{|v-(v)_{2\tau}-(Dv)_{\tau }y|^2}{\tau^2}\,dy,
\end{align*}
where we have used that $v$ and $Dv$ are bounded, $\lambda_ h\rightarrow 0$ and $|V(\xi)|\leq|\xi|$ for $|\xi|\leq 1$.\\
In view of this aim it is enough to prove that
\begin{equation}
I:= \lim_{h\to \infty}\frac{1}{\lambda_h^2}\medint_{B_{2\tau}}\bigg|V\bigg(\frac{\lambda_h\bigl((v_h-v) -(v_h-v)_{2\tau}-(Dv_h-Dv)_{\tau}\,y\bigr)}{2\tau}\bigg)\bigg|^2\,dy=0.   
\end{equation}
In the sequel  $\sigma$ will denote the exponent given in the Sobolev-Poincar\'e type inequality of the Theorem \ref{SoPoV}.
We can assume that the higher integrability exponent $\delta$ given in the step $5$ is such that $\delta<\sigma$.\\
\indent Let us choose $\theta\in(0,1)$ such that $2\theta+\frac{1-\theta}{1+\sigma}=1$. Applying H\"older's inequality, it holds that
\begin{align*}
&0\leq I\leq\limsup_{h\to \infty}\frac{1}{\lambda_h^2}\biggl(\medint_{B_{2\tau}}\bigg|V\bigg(\frac{\lambda_h\bigl((v_h-v) -(v_h-v)_{2\tau}-(Dv_h-Dv)_{\tau}\,y\bigr)}{2\tau}\bigg)\bigg|\,dy\biggr)^{2\theta}\\
& \cdot \biggl(\medint_{B_{2\tau}}\bigg|V\bigg(\frac{\lambda_h\bigl((v_h-v) -(v_h-v)_{2\tau}-(Dv_h-Dv)_{\tau}\,y\bigr)}{2\tau}\bigg)\bigg|^{2(1+\sigma)} dy\biggr)^{\frac{1-\theta}{1+\sigma}}.
\end{align*}
Using the fact that $|V(\xi)|\leq |\xi|$ and \emph{(iii)} of Lemma \ref{PropV}, for the first factor in the previous product, and using also Theorem \ref{SoPoV} applied to $(v_h-v) -(v_h-v)_{2\tau}-(Dv_h-Dv)_{\tau}\,y$, we deduce 
\begin{align*}
&0\leq I\leq\limsup_{h\to \infty}\frac{c}{\lambda_h^2}\bigg(\lambda_h\medint_{B_{2\tau}}\bigg(\bigg|\frac{v_h-v}{\tau}\bigg|+\bigg|\frac{(Dv_h-Dv)_{\tau}}{\tau}\bigg|\bigg)\,dy\bigg)^{2\theta}\\
& \times \bigg(\medint_{B_{6\tau}}\big|V\big(\lambda_h(Dv_h-Dv)\big) \big|^{\alpha} dy\bigg)^{\frac{2(1-\theta)}{\alpha}}.
\end{align*}
In the last term we can increase choosing $\alpha=2$ and accordingly, observing that
\begin{equation*}
\int_{B_1}\big|V\big(\lambda_h Dv_h\big) \big|^{2} dy\leq c(n)\lambda_h^2,
\end{equation*}
we conclude that
\begin{align*}
&0\leq I\leq\lim_{h\to \infty}\frac{c}{\lambda_h^2}\lambda_h^{2\theta}\bigg(\medint_{B_{2\tau}}\bigg(\bigg|\frac{v_h-v}{\tau}\bigg|+\bigg|\frac{(Dv_h-Dv)_{\tau}}{\tau}\bigg|\bigg)\,dy\bigg)^{2\theta}\cdot c\lambda_h^{2(1-\theta)}\\
& =\lim_{h\to \infty}C\biggl(\medint_{B_{2\tau}}\big(|{v_h-v}|+|{(Dv_h-Dv)_{\tau}}|\big)\,dy\biggr)^{2\theta}=0.
\end{align*}
By virtue of \eqref{24bis}, \eqref{25quater},  \eqref{25ter}, \eqref{eqq2}, the   Poincar\'e-Wirtinger inequality and \eqref{eqq2}, we get
\begin{align*}
\limsup_{h\to \infty}\frac{U_*(x_h,\tau r_h)}{\lambda^2_h}
&\leq c(n,p,\ell_1,\ell_2,L_2,\Lambda,M) \bigg\{\medint_{B_{2\tau}}\frac{|v-(v)_{2\tau}-(Dv)_{\tau }y|^2}{\tau^2}\,dy+\tau\bigg\}\\
& \leq c(n,p,\ell_1,\ell_2,L_2,\Lambda,M)\bigg\{\medint_{B_{2\tau}}|Dv-(Dv)_{\tau}|^2\,dy+\tau\bigg\}\\
& \leq c(n,N,p,\ell_1,\ell_2,L_1,L_2,\Lambda,M)\big[\tau^2+\tau\big]\leq C(n,N,p,\ell_1,\ell_2,L_1,L_2,\Lambda,M)\tau.
\end{align*}
The contradiction follows,  by choosing $C_*$ such that $  C_*>C$, since, by  \eqref{decay4},
\begin{equation}
\liminf_h\frac{U_*(x_h,\tau r_h)}{\lambda^2_h}\ge C_*\tau.
\end{equation}
\end{proof}
If assumption \eqref{H} is not taken into account, it is still possible to establish a decay result for the excess, analogous to the previous one. However, this requires employing a modified "hybrid" excess, defined as:
\begin{equation}\label{excess2}U_{**}(x_0,r):=U(x_0,r)+ \left(\frac{P(E, B_r(x_0))}{r^{n-1}}\right)^{\frac{\delta}{1+\delta}}+r^\beta,
\end{equation}
where $U(x_0,r)$ is defined in \eqref{excess}, $\delta$ is the higher integrability exponent given in the Step 5 of Proposition \ref{decay1}  and $0<\beta<\frac{\delta}{1+\delta}$. The following result still holds true.
\medskip
\begin{Prop}\label{decay1a}  Let $(u,E)$ be a local minimizer of $\mathcal{I}$ in \eqref{intro0} under the assumptions \eqref{F1p}, \eqref{F2p}, \eqref{G1p} and \eqref{G2p}. For every $M>0$ and $0<\tau <\frac{1}{4}$, there exist two positive constants $\varepsilon_0=\varepsilon_0(\tau,M)$ and $c_{**}=c_{**}(n,p,\ell_1,\ell_2,L_1,L_2,\Lambda,\delta,M)$ for which, whenever $B_r(x_0)\Subset\Om$ verifies
\begin{equation*}
    |(Du)_{x_0,r}|\leq M\quad\mathrm{and}\quad  U_{**}(x_{0}, r)\leq \epsilon_0,
\end{equation*}
 then
\begin{equation}\label{decay3a}
U_{**}(x_{0}, \tau r)\leq c_{**}\, \tau^\beta\, U_{**}(x_{0}, r).
 \end{equation}
 \end{Prop}
 In order to avoid unnecessary repetition we do not include the proof here, as it is almost identical to the proof of the Proposition \ref{decay1}, with the obvious adjustments, see \cite{CEL}.
\bigskip
\section{Proof of the Main Theorem}
Here we give the proof of Theorem \ref{main} through a suitable iteration procedure. It is easy to show the validity of the following lemma by arguing exactly in the same way as in \cite[Lemma 6.1]{CFP1}.  
\begin{Lem}\label{ite} 
Let $(u,E)$ be a minimizer of the functional $\mathcal{I}$ and let $c_*$  the constant introduced in Proposition \ref{decay1}. For every $\alpha\in (0,1)$ and $M>0$ there exists $\vartheta_0=\vartheta_0(c_*,\alpha)<1$ such that for $\vartheta\in (0,\vartheta_0)$ there exists a positive constant $\varepsilon_1=\varepsilon_1(n,p,\ell_1,\ell_2,L_1,L_2,M,\vartheta)$ such that, if $B_r(x_0)\Subset \Omega$,
\begin{equation*}
    |Du|_{x_0,r}<M\quad \text{and}\quad U_*(x_0,r)<\varepsilon_1,
\end{equation*}
 then
\begin{equation}
\label{indu2}
|D u|_{x_0,\vartheta^h r}<2M\quad\text{and}\quad
U_*(x_0,\vartheta^{h}r) \leq\vartheta^{h\alpha} U_*(x_0,r),\quad \forall h\in \mathbb{N}_0.
\end{equation}

\end{Lem}
\begin{proof}
Let $M>0$, $\alpha\in (0,1)$ and $\vartheta\in (0,\vartheta_0)$, where $\vartheta_0<1$. Let $\varepsilon_1<\varepsilon_0$, where $\varepsilon_0$ is the constant appearing in Proposition \ref{decay1}. We first prove by induction that
\begin{equation}\label{indu1}
    |D u|_{x_0,\vartheta^h r}<2M, \quad\forall h\in\N_0.
\end{equation}
If $h=0$, the statement holds. Assuming that \eqref{indu1} holds for $h>0$, applying properties \emph{(i)} and \emph{(vi)} of Lemma \ref{PropV}, we compute: 
\begin{align}\label{it1}
|Du|_{x_0,\vartheta^{h+1} r}
& \leq  |Du|_{x_0,r}+\sum_{j=1}^{h+1}||Du|_{x_0,\vartheta^{j} r}-|Du|_{x_0,\vartheta^{j-1}r}|\\
& \leq M+\sum_{j=1}^{h+1} \medint_{B_{\vartheta^j r}}|Du-(Du)_{x_0,\vartheta^{j-1}r}|\,dx\\
& \leq M +\vartheta^{-n}\sum_{j=1}^{h+1}\Bigg[ \frac{1}{|B_{\vartheta^{j-1} r}|}\int_{B_{\vartheta^{j-1} r}\cap\{|Du-(Du)_{x_0,\vartheta^{j-1}r}|\leq 1\}}\!\!\!\!\!\!\!\!\!\!\!\!\!\!\!\!\!\!\!\!\!\!\!\!\!\!\!\!|Du-(Du)_{x_0,\vartheta^{j-1}r}|\,dx\\
& + \frac{1}{|B_{\vartheta^{j-1} r}|}\int_{B_{\vartheta^{j-1} r}\cap\{|Du-(Du)_{x_0,\vartheta^{j-1}r}|> 1\}}|Du-(Du)_{x_0,\vartheta^{j-1}r}|\,dx\Bigg]\\
& \leq M+\vartheta^{-n}\sum_{j=1}^{h+1} \Bigg[\bigg(\medint_{B_{\vartheta^{j-1} r}}|V(Du-(Du)_{x_0,\vartheta^{j-1}r})|^2\,dx\bigg)^{\frac{1}{2}}\\
& +\bigg(\medint_{B_{\vartheta^{j-1} r}}|V(Du-(Du)_{x_0,\vartheta^{j-1}r})|^2\,dx\bigg)^{\frac{1}{p}}\Bigg]\\
& \leq M+c(p,M)\vartheta^{-n}\sum_{j=1}^{h+1}\big[U_*(x_0,\vartheta^{j-1}r)^{\frac{1}{2}}+U_*(x_0,\vartheta^{j-1}r)^{\frac{1}{p}}\big]\\
& \leq M+c(p,c_*,M)\varepsilon_1^{\frac{1}{2}}\vartheta^{-n}\sum_{j=1}^{h+1} \vartheta^{\frac{j-1}{2}}\leq M+c(p,c_*,M)\varepsilon_1^{\frac{1}{2}}\frac{\vartheta^{-n}}{1-\vartheta^{\frac{1}{2}}}\leq 2M,
\end{align}
where we have chosen $\varepsilon_1=\varepsilon_1(p,c_*,M,\vartheta)>0$ sufficiently small. Now we prove the second inequality in \eqref{indu2}. The statement is obvious for $h=0$. If $h>0$ and \eqref{indu2} holds, we have that
\begin{equation}\label{it2}U_*(x_0,\vartheta^{h}r)\leq \vartheta^{h\alpha} U_*(x_0,r)<\varepsilon_1,\end{equation}
by our choice of $\vartheta$  and $\varepsilon_1$. Thus thanks to \eqref{indu1} we can apply Proposition \ref{decay1} with $\vartheta^h r$ in place of $r$  to deduce that
\begin{equation*}
    U_*(x_0,\vartheta^{h+1}r)\leq \vartheta^\alpha U_*(x_0,\vartheta ^h r)\leq \vartheta^{(h+1)\alpha}U_*(x_0, r),
\end{equation*}
where we have chosen  $\vartheta_0=\vartheta_0(c_*,\alpha)$ sufficiently small and we have used \eqref{it2}. Therefore, the second inequality in \eqref{indu2} is also true for every $k\in \mathbb{N}$.
\end{proof}
Analogously, it is possible to prove an iteration lemma for $U_{**}$.
\begin{Lem}\label{ite1} Let $(u,E)$ be a minimizer of the functional $\mathcal{I}$ and let $\beta$ be the exponent of Proposition \ref{decay1a}. For every $M>0$ and $\vartheta\in (0,\vartheta_0)$, with $\vartheta_0<\min\left\{ c_{**},\frac{1}{4}\right\}$, there exist $\varepsilon_1>0$ and $R>0$ such that, if   $r<R$ and $x_0\in\Omega$ satisfy
\begin{equation*}
    B_r(x_0)\Subset \Omega,\quad |Du|_{x_0,r}<M\quad and\quad U_{**}(x_0,r)<\varepsilon_1,
\end{equation*}
where $c_{**}$ is the constant introduced in Proposition \ref{decay1a}, then
\begin{equation}
|D u|_{x_0,\vartheta^h r}<2M\quad\text{and}\quad    U_{**}(x_0,\vartheta^{k}r)\leq \vartheta^{k\beta} U_{**}(x_0,r), \quad\forall k\in\N.
\end{equation}

\end{Lem}
\smallskip
\begin{proof}[Proof of Theorem \ref{main}]

\noindent We consider the set 
\begin{equation}
\Omega_1:=\bigg\{x\in \Omega:\,\, \limsup_{\rho\to 0}|(Du)_{x,\rho}|<\infty \,\, \mathrm{and}\,\, \limsup_{\rho\to 0} U_*(x,\rho)=0\bigg\}
\end{equation}
and let $x_0\in \Omega_1$. For every $M>0$ and for $\varepsilon_1$ determined in Lemma \ref{ite} there exists a radius $ R_{M,\varepsilon_1}>0$ such that
\begin{equation} |Du|_{x_0,r}<M\quad \text{and}\quad U_{*}(x_0,r)<\varepsilon_1,
\end{equation}
for every $0<r<R_{M,\varepsilon_1}$.
Let $0<\rho<\vartheta r<R$ and $h\in\mathbb{N}$ be such that $\vartheta^{h+1}r<\rho<\vartheta^h r$, where $\vartheta=\frac{\vartheta_0}{2}$ and $\vartheta_0$ is the same constant appearing in Lemma \ref{ite}. 
By  Lemma \ref{ite}, we obtain
\begin{equation*}
    |D u|_{x_0,\rho}\leq \frac{1}{\vartheta^n}|D u|_{x_0,\vartheta^h r}\leq c(M,c_*,\alpha).
\end{equation*}
Using \emph{(iv)} of Lemma \ref{PropV} and reasoning as in the proof of Lemma \ref{ite}, we estimate
\begin{align}
|V((Du)_{x_0,\vartheta^h r})-V((Du)_{x_0,\rho})|^2
& \leq c(n,p)|(Du)_{x_0,\vartheta^h r}-(Du)_{x_0,\rho}|^2\\
& \leq c(n,p,c_*,M)\vartheta_0^{-2n}\vartheta^{h\alpha}U_*(x_0,r).
\end{align}
Thus, taking the previous two inequalities into account, applying again Lemma \ref{ite}, we estimate
\begin{align}
& U_*(x_0,\rho)
\leq 2\medint_{B_\rho(x_0)}|Du-(Du)_{x_0,\vartheta^h r}|^2\,dx+2|(Du)_{x_0,\vartheta^h r}-(Du)_{x_0,\rho}|^2+ \frac{P(E,B_\rho(x_0))}{\rho^{n-1}}+\rho\\
& \leq c(n,p,M,c_*\vartheta_0)\bigg[\medint_{B_{\vartheta^h r}(x_0)}|Du-(Du)_{x_0,\vartheta^h r}|^2\,dx+\vartheta^{h\alpha}U_*(x_0,r)+\frac{P(E,B_{\vartheta^h r}(x_0))}{(\vartheta^h r)^{n-1}}+\vartheta^h r\bigg]\\
& \leq c(n,p,c_*,M,\vartheta_0) \big[U_*(x_0,\vartheta^h r)+\vartheta^{h\alpha} U_*(x_0, r)\big]\leq c(n,p,c_*,M,\vartheta_0)\left(\frac{\rho}{r}\right)^\alpha U_*(x_0, r).
\end{align}
\noindent The previous estimate implies that
\begin{equation}
    \label{eccesso}
   {U(x_0,\rho)}\leq C_*\left( \frac{\rho}{r}\right)^{\alpha}U_*(x_0, r),
\end{equation}
where $C_{*}=C_{*}(n,p,c_*,M,\vartheta_0)$.
Since $U_*(y,r)$ is continuous in $y$, we have that $U_*(y,r)<\varepsilon_1$ for every $y$ in a suitable neighborhood $I$ of $x_0$. Therefore, for every $y\in I$ we have that
\begin{equation*}
  U(y,\rho)\leq C_* \left(\frac{\rho}{r}\right)^\alpha U_*(y, r).  
\end{equation*}
The last inequality implies,
by the Campanato characterization of H\"older continuous  functions (see \cite[Theorem 2.9]{gi}), that $u$ is $C^{1,\alpha}$ in $I$ for every $0<\alpha<\frac{1}{2}$, and we can conclude that 
the set $\Omega_1$ is open and the function $u$ has  H\"older continuous  derivatives in $\Omega_1$.\\
\indent When the assumption \eqref{H} is not enforced, the proof goes exactly in the same way provided we use Lemma \ref{ite1} in place of Lemma \ref{ite}, with
 \begin{equation}
 \Omega_0:=\bigg\{x\in \Omega:\,\, \limsup_{\rho\to 0}|(Du)_{x_0,\rho}|<\infty \,\, \mathrm{and}\,\, \limsup_{\rho\to 0} U_{**}(x_0,\rho)=0\bigg\}.
 \end{equation}
\end{proof}

\noindent {\bf Acknowledgments:} The authors are members of the Gruppo Nazionale per l’Analisi Matematica, la Probabilità e le loro Applicazioni (GNAMPA) of the Istituto Nazionale di Alta Matematica (INdAM) and wish to acknowledge financial support from INdAM GNAMPA
Project 2024 “Regolarità per problemi a frontiera libera e disuguaglianze funzionali in
contesto finsleriano”.

\smallskip
\noindent

\noindent {\bf Conflict of interest:} The authors state no conflict of interest.

\smallskip
\noindent


\noindent
\author{Menita Carozza}, Dipartimento di Ingegneria, Universit\`{a} del Sannio, Corso Garibaldi,  Benevento 82100, Italy

\noindent
carozza@unisannio.it
\medskip

\noindent
\author{Luca Esposito},
Dipartimento di Matematica, Università degli Studi di Salerno, Via Giovanni Paolo II 132, Fisciano 84084, Italy\\
luesposi@unisa.it

\medskip
\noindent
\author{Lorenzo Lamberti},
Dipartimento di Matematica, Università degli Studi di Salerno, Via Giovanni Paolo II 132, Fisciano 84084, Italy\\
llamberti@unisa.it

\end{document}